\documentclass[11pt]{amsart}
\usepackage[margin=1.35in]{geometry}
\usepackage{amscd,amsmath,amsxtra,amsthm,amssymb,stmaryrd,xr,mathrsfs,mathtools,enumerate,commath, comment}
\usepackage{stmaryrd}
\usepackage{xcolor}
\usepackage{commath}
\usepackage{comment}
\usepackage{tikz-cd}
\usepackage{float}
\usepackage{longtable} 
\usepackage{pdflscape} 
\usepackage{booktabs}
\usepackage{hyperref}
\definecolor{vegasgold}{rgb}{0.77, 0.7, 0.35}
\definecolor{darkgoldenrod}{rgb}{0.72, 0.53, 0.04}
\definecolor{gold(metallic)}{rgb}{0.83, 0.69, 0.22}
\hypersetup{
 colorlinks=true,
 linkcolor=darkgoldenrod,
 filecolor=brown,      
 urlcolor=gold(metallic),
 citecolor=darkgoldenrod,
 pdftitle={On the Distribution of Iwasawa invariants associated to multigraphs},
 }

\usepackage[all,cmtip]{xy}

\usepackage{tikz}
\usetikzlibrary{shapes.geometric}
\usetikzlibrary{decorations.markings}
\tikzset{every loop/.style={min distance=10mm,looseness=10}}

\DeclareFontFamily{U}{wncy}{}
\DeclareFontShape{U}{wncy}{m}{n}{<->wncyr10}{}
\DeclareSymbolFont{mcy}{U}{wncy}{m}{n}
\DeclareMathSymbol{\Sh}{\mathord}{mcy}{"58}
\usepackage[T2A,T1]{fontenc}
\usepackage[OT2,T1]{fontenc}

\newtheorem{theorem}{Theorem}[section]
\newtheorem{lemma}[theorem]{Lemma}
\newtheorem{question}[theorem]{Question}

\newtheorem{proposition}[theorem]{Proposition}
\newtheorem{corollary}[theorem]{Corollary}
\newtheorem{definition}[theorem]{Definition}
\newtheorem{assumption}[theorem]{Assumption}

\numberwithin{equation}{section}

\theoremstyle{remark}
\newtheorem{remark}[theorem]{Remark}
\newtheorem{example}[theorem]{Example}

\newcommand{\FF}{\mathbb{F}}

\newcommand{\ord}{\mathrm{ord}}

\newcommand{\lb}{[\![}
\newcommand{\rb}{]\!]}

\newcommand{\Z}{\mathbb{Z}}

\newcommand{\Q}{\mathbb{Q}}
\newcommand{\F}{\mathbb{F}}

\newcommand{\Zl}{\mathbb{Z}_\ell}

\newcommand{\blue}[1]{\textcolor{blue}{#1}}

\newcommand{\PP}{\op{Prob}}
\newcommand{\Fl}{\FF_\ell}

\newcommand{\op}[1]{\operatorname{#1}}

\newcommand\mtx[4] { \left( {\begin{array}{cc}
 #1 & #2 \\
 #3 & #4 \\
 \end{array} } \right)}

\begin{document}
\title[On the distribution of Iwasawa invariants associated to multigraphs]{On the distribution of Iwasawa invariants associated to multigraphs}

\author[C.~Dion]{C\'{e}dric Dion}
\address{C\'{e}dric Dion\newline D\'epartement de Math\'ematiques et de Statistique\\
Universit\'e Laval, Pavillion Alexandre-Vachon\\
1045 Avenue de la M\'edecine\\
Qu\'ebec, QC\\
Canada G1V 0A6}
\email{cedric.dion.1@ulaval.ca}

\author[A.~Lei]{Antonio Lei}
\address{Antonio Lei\newline Department of Mathematics and Statistics\\University of Ottawa\\
150 Louis-Pasteur Pvt\\
Ottawa, ON\\
Canada K1N 6N5}
\email{antonio.lei@uottawa.ca}

\author[A.~Ray]{Anwesh Ray}
\address{Anwesh Ray\newline Department of Mathematics\\
University of British Columbia\\
Vancouver BC, Canada V6T 1Z2}
\email{anweshray@math.ubc.ca}

\author[D.~Valli\`{e}res]{Daniel Valli\`{e}res}
\address{Daniel Valli\`{e}res\newline Mathematics and Statistics Department, California State University, Chico, CA 95929 USA}
\email{dvallieres@csuchico.edu}

\begin{abstract}
Let $\ell$ be a prime number. The Iwasawa theory of multigraphs is the systematic study of growth patterns in the number of spanning trees in abelian $\ell$-towers of multigraphs. In this context, growth patterns are realized by certain analogues of Iwasawa invariants, which depend on the prime $\ell$ and the abelian $\ell$-tower of multigraphs. We formulate and study statistical questions about the behaviour of the Iwasawa $\mu$ and $\lambda$ invariants. 
\end{abstract}

\date{\today} 

\maketitle
\tableofcontents 
 
\section{Introduction}\addtocontents{toc}{\protect\setcounter{tocdepth}{1}}
\subsection*{Classical Iwasawa theory} Let $\ell$ be a prime number, $K$ be a number field and $\Z_\ell$ denote the ring of $\ell$-adic integers. A \emph{$\Z_\ell$-tower} over $K$ consists of a sequence of abelian extensions \[K=K_0\subset K_1\subset K_2\subset K_3\subset \dots \subset K_n\subset K_{n+1}\subset \dots \]
such that for all values of $n$, the Galois group $\op{Gal}(K_n/K)$ is isomorphic to $\Z/\ell^n\Z$. The union $K_\infty:=\bigcup_{n\geq 1} K_n$ is an infinite Galois extension of $K$ with Galois group $\op{Gal}(K_\infty/K)$ isomorphic to $\Z_\ell$. We refer to $K_\infty$ as a \emph{$\Z_\ell$-extension} of $K$. For $n\geq 0$, we denote by $\op{Cl}_{\ell}(K_n)$, the Sylow $\ell$-subgroup of the ideal class group of $K_n$, and by $h_\ell(K_n)$ the cardinality of $\op{Cl}_{\ell}(K_n)$. In his seminal work \cite{iwasawa1973zl}, Iwasawa studied growth patterns for $h_\ell(K_n)$, as $n\rightarrow \infty$. More specifically, it is shown in \emph{loc. cit.} that there are invariants $\mu$, $\lambda$ and $\nu$, which only depend on the prime $\ell$ and the extension $K_\infty$ (and independent of $n$), such that for $n\gg 0$,
\begin{equation}\label{eq:Iw}
    h_{\ell}(K_n)=\ell^{(\mu \ell^n+ \lambda \ell+\nu)}.
\end{equation} The quantities $\mu$ and $\lambda$ are non-negative integers, while the quantity $\nu$ is an integer, possibly negative.

\par Of special interest is the case when $K_\infty/K$ is the \emph{cyclotomic $\Z_\ell$-extension} of $K$, in which case we simply denote the Iwasawa invariants by $\mu_\ell(K)$, $\lambda_\ell(K)$ and $\nu_\ell(K)$. By the celebrated result of Ferrero and Washington \cite{ferrero1979iwasawa}, $\mu_\ell(K)$ is equal to $0$ for all abelian extensions $K/\Q$. Moreover, it is conjectured in \cite{Iw73} that $\mu_\ell(K)$ is $0$ for all number fields $K$ . The $\lambda$-invariant is more subtle, and it is natural to study statistical questions about its average behavior. 

Let $n>0$ be an integer, let $\mathcal{F}$ be a family of number fields of fixed degree $n$. We do not provide an explicit definition of what constitutes a \emph{family}, however it is natural to study number fields satisfying prescribed splitting conditions at a fixed set of primes. For instance, it is natural to consider the case when $\mathcal{F}$ is the family of imaginary quadratic fields in which $p$ is inert (or split). 
\par For a number field $K$, we let $|\Delta_K|$ be the absolute discriminant. Given a real number $x>0$, let $\mathcal{F}_{\leq x}$ be the subset of $\mathcal{F}$ consisting of number fields $K$ such that $|\Delta_K|\leq x$. Given $\mu_0\geq 0$ and $\lambda_0\geq 0$ and a real number $x>0$, we let $\mathcal{F}_{\leq x}(\mu_0, \lambda_0, \ell)$ consist of $K\in \mathcal{F}_{\leq x}$ such that $\mu_\ell(K)=\mu_0$ and $\lambda_\ell(K)=\lambda_0$. Note that both $\mathcal{F}_{\leq x}$ and $\mathcal{F}_{\leq x}(\mu_0, \lambda_0, \ell)$ are finite sets. The following question motivates the line of inquiry in this manuscript.
\begin{question}\label{Q}
Fix a prime number $\ell$. Can one compute the proportion of number fields $K$ in $\mathcal{F}$ such that $\mu=\mu_0$ and $\lambda=\lambda_0$. In greater detail, one is interested in computing the density defined as follows
\[\mathfrak{d}(\mu_0, \lambda_0, \ell)=\lim_{x\rightarrow \infty} \frac{\#\mathcal{F}_{\leq x} (\mu_0, \lambda_0, \ell)}{\# \mathcal{F}_{\leq x}},\]in cases when the above limit exists.
\end{question}
This can be regarded as an analogue of  similar questions on classical Iwasawa invariants of number fields. The heuristics for the distribution of these arithmetic Iwasawa invariants have been studied by Ellenberg, Jain and Venkatesh \cite{ellenberg2011modeling}. One may attempt to compute the limit above for various choices of $(\mu_0, \lambda_0)$. Of course, the $\mu=0$ conjecture would imply that $\mathfrak{d}(\mu_0, \lambda_0, \ell)=0$ if $\mu_0>0$. It turns out that the above question, though of considerable interest, leads to difficulties that are yet to be resolved systematically.

\subsection*{Iwasawa theory of multigraphs}
In the present article, we consider  graph-theoretic analogs of Iwasawa invariants, for which  questions similar to Question~\ref{Q} become tractable for certain families of multigraphs. What might be called the Iwasawa theory of multigraphs is a recent area of research. The reader may refer to \cite{Vallieres:2021, mcgownvallieresII, mcgownvallieresIII} where the approach via Artin-Ihara $L$-functions is introduced and to \cite{Gonet:2022} where the module theoretical approach was initiated.  Other papers on the subject where the reader can find more information are \cite{lei2022non}, \cite{DuBose/Vallieres:2022}, \cite{Kleine/Muller:2022}, \cite{Ray/Vallieres:2022}, and \cite{Kataoka:2023}. However, at present, there is no standard text on the subject.

A finite multigraph $X$ consists of a finite set of vertices joined by directed edges and an inversion map satisfying prescribed properties. For the precise definition of a multigraph, the reader is referred to Definition \ref{definition multigraph}. In this article, we allow there to be more than one edge joining any two vertices. We  also allow for loops based at a single vertex. It is assumed that there are only finitely many edges, and the multigraph is connected. Given a prime $\ell$, a \emph{connected $\Z_\ell$-tower}
over $X$ consists of a sequence of covers of connected multigraphs
\[X=X_0\leftarrow X_1\leftarrow X_2\leftarrow \dots \leftarrow X_n \leftarrow X_{n+1}\leftarrow \dots, \]
such that for every integer $n\geq 0$, the Galois group $\op{Gal}(X_n/X)$ is isomorphic to $\Z/\ell^n \Z$. Note that such a tower is not necessarily unique. The multigraphs $X_n$ can be constructed from a $\Zl$-valued voltage assignment on the edges of $X$. For further details, please refer to \S\ref{S:review} where this construction is summarized. We establish a sufficient condition for $X_n$ to be connected for all $n\ge0$ (see Theorem~\ref{connectedness_cond}). This particular criterion is potentially of independent interest. The inverse limit $\varprojlim_n \op{Gal}(X_n/X)$ is isomorphic to $\Z_\ell$, and thus a $\Zl$-tower of a multigraph is analogous to  a $\Z_\ell$-extension of a number field. The Picard group $\op{Pic}^0(X_n)$ has cardinality equal to the number of spanning trees in $X_n$, for more details, cf. \cite[Theorem 10.8]{Sunada:2013}. The cardinality of the Sylow $\ell$-subgroup of $\op{Pic}^0(X_n)$ is denoted by $\kappa_\ell(X_n)$. It is shown by McGown and the fourth named author (cf. \cite{Vallieres:2021,mcgownvallieresII,mcgownvallieresIII}) that there are invariants $\mu, \lambda\in \Z_{\geq 0}$ and $\nu\in \Z$ such that 
\begin{equation}\label{eq:Iw-graphs}
    \kappa_\ell(X_n)=\ell^{(\mu\ell^n +\lambda \ell+\nu)}
\end{equation} for $n\gg 0$. This establishes an analog of Iwasawa's formula \eqref{eq:Iw} in the setting of multigraphs. The formula \eqref{eq:Iw-graphs} above was also independently obtained by Gonet, see \cite{Gonet:2021a, Gonet:2022}. We prove  that when $\ell$ is odd, the invariant  $\lambda$ is in fact always an odd integer (see Theorem~\ref{lambda is odd}), which again may be of independent interest.

\subsection*{Main results} The main goal of the present paper is to study the distribution of Iwasawa $\mu$ and $\lambda$-invariants for multigraphs that feature in \eqref{eq:Iw-graphs}. We study two related questions in this context. First, we fix a multigraph $X$ and allow the $\Z_\ell$-tower over $X$ to vary in a natural way and study the variation of Iwasawa invariants. These investigations are carried out for multigraphs with one or two vertices as well as certain families of $\Zl$-towers of complete graphs. We expect that the methods developed in the present paper can be generalized to many more families of multigraphs. However, the calculations may become cumbersome for general multigraphs. The second theme is more closely related to arithmetic statistics. Here, we consider all multigraphs with a given number of vertices, and their connected $\Z_\ell$-towers and study the variation of $\mu$ and $\lambda$-invariants for this entire family. This question can be formulated for multigraphs with a fixed number of vertices $N$ for any number $N$. However, we find that the computations are significantly involved even in the case when $N=1$. We refer to question \eqref{bouquet question 2} for the precise formulation of the problem in this setting. Below is a more detailed description of the results we prove.
\setlength{\leftmargini}{1cm}

\begin{enumerate}[(A)]
\item In \S\ref{S:bouquet}, we give a detailed analysis of the setting where $X$ is a bouquet (a multigraph with a single vertex). In particular:
\begin{enumerate}[({A}.1)]
    \item We first establish several sufficient and/or necessary conditions for $\mu=0$. This allows us to prove that $\mu$ is always zero and give an upper bound on $\lambda$ under the hypothesis that $\ell$ is larger than the number of loops in a bouquet   (see Corollary~\ref{cor:posmu} and Theorem~\ref{thm:small-t} for the precise statements). Without this hypothesis,  we show that both $\mu$ and $\lambda$  can be arbitrarily large (see Lemma~\ref{arb_large}), which mirrors similar results in the setting of  Iwasawa theory of number fields and elliptic curves \cite{Iw73,matsuno,kim09,kr22}. 
    \item Let $t$ be the number of loops in a bouquet. We study the distribution of $\mu$ and $\lambda$ as the $\Zl$-tower varies. As we shall see in \S\ref{S:review}, where we review the Iwasawa theory of multigraphs, the connected $\Zl$-towers of a bouquet with $t$ loops are in one-one correspondence  with $\Zl^t\setminus (\ell\Zl)^t$. In particular, the Haar measure on $\Zl$ gives a measure on the set of all connected $\Zl$-towers of a given bouquet. Under this measure, we are able to obtain an upper bound on the probability for $\mu>0$, and for $\mu=0$ with a fixed choice of $\lambda$ which is sufficiently small (see Theorems~\ref{thm:mu>0} and \ref{thm:mu0lambda small}). For $\mu=0$ and $\lambda=1$, we are able to prove a precise formula  (see Theorem~\ref{thm:mu0lambda1}). We also investigate the cases where $t=2,3$ separately in \S\ref{S:smallt}. 
    \item Specializing to the case where the $\Zl$-towers correspond to the vectors $\alpha=(\alpha_1,\ldots,\alpha_t)\in \Z^t\setminus(\ell\Z)^t$, we  show in Theorem~\ref{section 3 main thm} that as both $t$ and $\alpha$ vary, the probability of  $\mu=0$ and $\lambda=1$ is equal to $1-\ell^{-1}$, mirroring a similar result obtained by the first and third named authors in the context of 2-bridge links (see \cite[Theorem~5.4]{dionray}). 
\end{enumerate}
\item In \S\ref{S:2-vertices}, we study  multigraphs with two vertices. In particular, we determine a precise formula for the probability of a $\Zl$-tower of a given multigraph with two vertices to have $\mu=0$ and $\lambda=1$ via the theory of quadratic forms over finite fields. See in particular Theorem~\ref{main thm 2 vertices}. 
\item Finally, we study $\Zl$-towers of complete graphs. In particular,   Theorem~\ref{thm:complete}  gives a precise description of the distribution of $\mu$ and $\lambda$ for certain families including the so-called cyclic voltage covers of the complete graph $K_u$ with $u$ vertices (which were initially studied by Gonet \cite{Gonet:2021}) as $u$ varies.
\end{enumerate}

\subsection*{Outlook} The present paper initiates the study of  distributions of Iwasawa invariants of multigraphs and answers several statistical questions in some special cases. Our main results rely crucially on studying the power series in $\Zl\lb T\rb$ that governs the growth formula \eqref{eq:Iw-graphs}. For example, in order to study when $\mu=0$ and $\lambda=1$ occurs, we are faced with the task of analyzing whether the leading term of the aforementioned power series is divisible by $\ell$. In the cases being treated, we are able to reformulate this question in terms of explicit quadratic forms over $\Fl$ and employ results on finite fields  (e.g. the Chevalley--Warning theorem).  On developing further results on systems of equations over finite fields in  several variables, one may study more than one coefficient of a power series in $\Zl\lb T\rb$ simultaneously. This  should then allow us to extend our results on $\mu=0$ and $\lambda=1$ to more general combinations of Iwasawa invariants. As mentioned above, carrying out similar calculations beyond the cases we have studied so far can potentially be cumbersome. It would  be desirable to develop a systematic approach to study the questions inaugurated in the present paper for more general families of graphs and multigraphs. For example, it would be very interesting to extend our methods to study distributions of Iwasawa invariants that arise from Erd\"os--R\'enyi random graphs. We plan to study some of these questions in the future.

\subsection*{Acknowledgement}CD's research is supported by the NSERC Canada Graduate Scholarships – Doctoral program. AL's research is supported by the NSERC Discovery Grants Program RGPIN-2020-04259 and RGPAS-2020-00096. DV would like to thank Kevin McGown for several stimulating discussions. We thank the anonymous referee for valuable comments and suggestions.

\section{Preliminaries}
\par This section is dedicated to the discussion of preliminary notions in this paper. 
\addtocontents{toc}{\protect\setcounter{tocdepth}{2}}
\subsection{Basic properties of multigraphs}
Fix a prime number $\ell$. We introduce basic notions in the theory of multigraphs and Iwasawa theory of $\Z_\ell$-towers. For further details, the reader may refer to \cite[\S4]{mcgownvallieresIII} and \cite[\S2]{lei2022non}. 
A finite undirected graph $X$ consists of a finite set $V_X$ of vertices $\{v_1, \dots, v_n\}$, and a finite number of edges $E_X$ joining the vertices. For such a graph, there is at most one edge which joins $v_i$ and $v_j$. Thus, in this context, we may view $E_X$ as a subset of $V_X\times V_X$. More specifically, if there is an edge joining $v_i$ to $v_j$, then $(v_i, v_j)$ belongs to $E_X$. Since the graph $X$ in question is undirected, we have that $(v_i, v_j)\in E_X$ if and only if $(v_j, v_i)\in E_X$. An edge from $v_i$ to $v_i$ is referred to as a \emph{loop}. The graph $X$ is fully described by the pair $(V_X, E_X)$. On the other hand, a directed graph is a graph for which there is at most one directed edge from $v_i$ to $v_j$. In this case, the set of edges is denoted $E_X^+$, and is a subset of $V_X\times V_X$. More specifically, there is a directed edge joining $v_i$ to $v_j$ precisely when $(v_i, v_j)\in E_X^+$. A \emph{multigraph} is a graph with possibly multiple directed edges joining any two vertices, subject to some further conditions. In this case, the set $E_X^+$ is not a subset of $V_X\times V_X$, but simply admits a map to $V_X\times V_X$. The notion of a multigraph is more general than that of a graph. We give the precise definition below.
\begin{definition}\label{definition multigraph}A finite multigraph $X$ is the data consisting of a finite set of vertices $V_X$ and a finite set of \emph{directed edges} $E_{X}^+$ along with functions
\[\begin{split}& i:E_{X}^{+} \to V_{X} \times V_{X},\\
&\iota:E_{X}^{+} \to E_{X}^{+}.
\end{split}\]
Here, $i$ and $\iota$ are the \emph{incidence} and \emph{inversion} functions, respectively, and are subject to the following conditions
\begin{enumerate}
\item $\iota^{2} = {\rm id}_{E_{X}^{+}}$,
\item $\iota(e) \neq e$ for all $e \in E_{X}^{+}$,
\item $i(\iota(e)) = \tau(i(e))$ for all $e \in E_{X}^{+}$,
\end{enumerate}
where $\tau(x,y) := (y,x)$.
\end{definition}
Given a directed edge $e\in E_X^+$, the incidence function maps $e$ to $i(e)=(v,v')$. The edge starts at $v$ and ends at the vertex $v'$. We define functions $o,t: E_X^+\rightarrow V_X$, such that $i(e)=(o(e), t(e))$. In other words, $o(e)$ (resp. $t(e)$) is the projection of $i(e)$ to the first (resp. second) coordinate. We refer to the function $o$ (resp. $t$) as the source (resp. target) function. The inversion map $\iota$ satisfies the property that $i(\iota(e))=(t(e), o(e))$. We shall occasionally denote the inverse $\iota(e)$ by $\bar{e}$. These two conventions will be used  interchangeably throughout.

\par The set of directed edges $e$ starting at a vertex $v$ and ending at a vertex $v'$ are the edges which satisfy the relation $i(e)=(v,v')$. We note that there may be any number of directed edges joining two vertices $v$ and $v'$, thus the graph $X$ is referred to as a multigraph. We allow for there to be edges to start and end at the same vertex $v$, and such edges are referred to as \emph{loops} based at $v$. We note that for distinct vertices $v$ and $v'$, there are as many directed edges starting at $v$ and ending at $v'$ as there are from $v'$ to $v$. Moreover, given a vertex $v$, there are an even number of directed loops starting and ending at $v$.

\par \textbf{Example:} Let us illustrate the above definitions through an example. The set of vertices is $V_X=\{v_1, v_2\}$. There are a total of $6$ directed edges $E_X^+=\{a,b, c, \bar{a}, \bar{b}, \bar{c}\}$, as depicted in the figure below.




\begin{figure}[h]
\begin{center}
\begin{tikzpicture}

\draw[fill=black] (0,0) circle (1pt) node[below]{$v_{1}$};
\draw[fill=black] (6,0) circle (1pt) node[below]{$v_{2}$};

\draw (0,0) edge [decoration={markings, mark= at position 0.44 with {\arrow[xscale=-1]{stealth}}},preaction={decorate},bend left=20, "$\bar{c}$" ] (6,0);
\draw (0,0) edge [decoration={markings, mark= at position 0.55 with {\arrow{stealth}}},preaction={decorate},bend left=60, "$b$"] (6,0);

\draw (0,0) edge [decoration={markings, mark= at position 0.58 with {\arrow{stealth}}},preaction={decorate}, bend right=20, , "$c$"] (6,0);
\draw (0,0) edge [decoration={markings, mark= at position 0.47 with {\arrow[xscale=-1]{stealth}}},preaction={decorate}, bend right=60, "$\bar{b}$"] (6,0);

\draw (0,0) edge [decoration={markings, mark= at position 0.28 with {\arrow{stealth}}},preaction = {decorate}, loop left, in = 155, out = 205,min distance=16mm, "$\bar{a}$"] (0,0);
\draw (0,0) edge [decoration={markings, mark= at position 0.28 with {\arrow[xscale=-1]{stealth}}},preaction = {decorate}, loop left, in = 140, out = 220,min distance=32mm, "$a$"] (0,0);
\end{tikzpicture}   
 
\end{center}
\end{figure}
With respect to notation above, $i(a)=(v_1, v_1)$, $i(b)=i(c)=(v_1, v_2)$.

\subsection{Iwasawa theory of multigraphs}\label{S:review}
\par In this section, we discuss the Iwasawa theory of multigraphs. The basic object of study is a $\Z_\ell$-tower of multigraphs. Let $X$ be a multigraph. Recall from the previous section that the data associated with $X$ consists of a finite set of vertices $V_X=\{v_1, \dots, v_n\}$ and finitely many directed edges $E_X^+$, along with an incidence map $i:E_X^+\rightarrow V_X\times V_X$, and an inversion map $\iota: E_X^+\rightarrow E_X^+$. These maps are subject to further conditions. Before introducing $\Z_\ell$-towers of multigraphs, we first discuss the notion of a Galois cover of multigraphs. For $v\in V_X$, let 
\[E_{X, v}:=\{e\in E_X^+\mid o(e)=v\}.\]
Let $X$ and $Y$ be two multigraphs.  A morphism of graphs $f:Y \rightarrow X$ consists of two functions $f_{V}:V_{Y} \rightarrow V_{X}$ and $f_{E}:E_{Y}^+ \rightarrow E_{X}^+$ satisfying
\begin{enumerate}
\item $f_{V}(o(e)) = o(f_{E}(e))$,
\item $f_{V}(t(e)) = t(f_{E}(e))$,
\item $\overline{f_{E}(e)} = f_{E}(\bar{e})$, \label{trois}
\end{enumerate}
for all $e \in E_{X}^+$.  We will often write $f$ for both $f_{V}$ and $f_{E}$.
\begin{definition} \label{cover}
Let $X$ and $Y$ be two graphs and let $f:Y\rightarrow X$ be a morphism of graphs. If $f$ satisfies the following conditions
\begin{enumerate}
\item $f:V_{Y} \rightarrow V_{X}$ is surjective,
\item for all $w \in V_{Y}$, the restriction $f|_{E_{Y,w}}$ induces a bijection 
$$f|_{E_{Y,w}}:E_{Y,w} \stackrel{\approx}{\rightarrow}{E_{X,f(w)}}, $$
\end{enumerate}
then $f$ is said to be a cover.
\end{definition}

\par If $Y/X$ is a cover and $Y$ is connected, then there exists a positive integer $d$ such that $|f^{-1}(v)| = d$ for all $v \in V_{X}$.  The integer $d$ is called the degree of $Y/X$ and will be denoted by $[Y:X]$.  If $f:Y \rightarrow X$ is a cover, then we let
$${\rm Aut}_{f}(Y/X) = \{ \sigma \in {\rm Aut}(Y) \, : \, f \circ \sigma = f\}. $$
\begin{definition}
The cover $f:Y \rightarrow X$ is called Galois if the following two conditions are satisfied.
\begin{enumerate}
\item The graphs $X$ and $Y$ are connected.
\item The group ${\rm Aut}_{f}(Y/X)$ acts transitively on the fiber $f^{-1}(v)$ for all $v \in V_{X}$.
\end{enumerate}
\end{definition}
If $Y/X$ is a Galois cover, then $[Y:X] = |{\rm Aut}_{f}(Y/X)|$, and we write ${\rm Gal}_{f}(Y/X)$, or ${\rm Gal}(Y/X)$ if $f$ is understood, instead of ${\rm Aut}_{f}(Y/X)$.  One has the usual Galois correspondence between subgroups of ${\rm Gal}(Y/X)$ and equivalence classes of intermediate covers of $Y/X$.
\par Galois covers of a multigraph can be constructed from certain combinatorial data associated to a multigraph, as we now explain. Define an equivalence relation on $E_X^+$, by setting $e\sim e'$ if $e=e'$ or $e=\bar{e}'$. The set of equivalence classes $E_X:=E_X^+/\iota$ is the set of \emph{undirected edges}, and let $\pi:E_X^+\to E_{X}$ be the $2:1$ natural quotient map, mapping a directed edge to the associated undirected edge. Choose a section $\gamma:E_{X} \to E_{X}^+$ to $\pi$, and set $S:=\gamma(E_{X})$. Our main reference for voltage assignments and derived multigraphs is \cite{Gross/Tucker:2001}, but we rephrase everything in the formalism used in \cite{Sunada:2013}. 
\begin{definition}
   Let $G$ be a finite abelian group.  With respect to notation above, a \emph{voltage assignment} is a function $\alpha:S \to G$.
\end{definition} We extend $\alpha$ to all of $E_{X}^{+}$ by setting $\alpha(\bar{s}) = \alpha(s)^{-1}$.  To this data, one can associate a multigraph $X(G,S,\alpha)$ as follows.  The set of vertices is $V = V_{X} \times G$ and the set of directed edges is $E^{+} = E_{X}^{+} \times G$.  The directed edge $(e,\sigma)$ connects $(o(e),\sigma)$ to $(t(e),\sigma \cdot \alpha(e))$ (where we recall that $o(e)$ and $t(e)$ denote the \textit{origin} and the \textit{target} of the edge $e$, respectively), and the inversion map is given by 
$$\overline{(e,\sigma)} = (\bar{e},\sigma \cdot \alpha(e)).$$
If $G_{1}$ is another finite abelian group and $f:G \to G_{1}$ is a group morphism, then one obtains a morphism of multigraphs $f_{*}:X(G,S,\alpha) \to X(G_{1},S,f \circ \alpha)$ that is defined via
$$f_{*}(v,\sigma) = (v,f(\sigma)) \text{ and } f_{*}(e,\sigma) = (e,f(\sigma)). $$
If both $X(G,S,\alpha)$ and $X(G,S,f \circ \alpha)$ are connected and $f$ is surjective, then $f_{*}$ is a cover in the sense of Definition \ref{cover} and in fact is a Galois cover with group of covering transformations isomorphic to ${\rm ker}(f)$.  In particular, if $f:G \to \{ 1\}$ is the group morphism into the trivial group and both $X$ and $X(G,S,\alpha)$ are connected, then we get a Galois cover $f_{*}:X(G,S,\alpha) \to X$ with group of covering transformation isomorphic to $G$.
\par We now recall the notion of a $\Z_\ell$-tower of multigraphs.
\begin{definition} \label{tower}
Let $\ell$ be a rational prime, and let $X$ be a connected graph.  A $\mathbb{Z}_{\ell}$-tower above $X$ consists of a sequence of covers of connected graphs
$$X = X_{0} \leftarrow X_{1} \leftarrow \ldots \leftarrow X_{n} \leftarrow \ldots $$
such that for all positive integer $n$, the cover $X_{n}/X$ obtained from composing the covers $X_{n} \rightarrow X_{n-1} \rightarrow \ldots \rightarrow X_{1} \rightarrow X$ is Galois with Galois group ${\rm Gal}(X_{n}/X)$ isomorphic to $\mathbb{Z}/\ell^{n}\mathbb{Z}$.
\end{definition}
\par We specialize the above construction to the case when $G=\Z_\ell$.  
Set $\alpha_{/n}$ to denote the mod-$\ell^n$ reduction of $\alpha$. We obtain a $\mathbb{Z}_{\ell}$-tower above $X$:
$$X  \leftarrow X(\mathbb{Z}/\ell\mathbb{Z},S,\alpha_{/1}) \leftarrow X(\mathbb{Z}/\ell^{2}\mathbb{Z},S,\alpha_{/2}) \leftarrow \ldots \leftarrow X(\mathbb{Z}/\ell^{k}\mathbb{Z},S,\alpha_{/k}) \leftarrow \ldots.$$  
Set $t:=|S|$ and enumerate $S=\{s_1,\dots , s_t\}$, and let $\alpha_i$ denote $\alpha(s_i)$. We may view $\alpha$ as a vector $\alpha=(\alpha_1, \alpha_2,\dots, \alpha_t)\in \Z_\ell^t$. The following assumption on the multigraphs $X(\Z/\ell^n \Z, S, \alpha_{/n})$ is necessary to relate the growth patterns in Picard groups to Iwasawa invariants. In section \ref{connectedness assumption section}, we introduce conditions under which the following assumption is satisfied. 
\begin{assumption}\label{main assumption}
Let $\ell$ be a prime (with $\ell=2$ allowed) and let $\alpha: S\rightarrow \Z_\ell$ be a function. We assume that the derived multigraphs $X(\Z/\ell^n \mathbb{Z}, S, \alpha_{/n})$ are connected for all $n\geq 0$.
\end{assumption}
\textbf{Examples:} Let us illustrate the above construction of a $\Z_\ell$-tower through some explicit examples.
\begin{enumerate}
\item Let $X$ be a bouquet with two loops, and choose an orientation $S =\{s_{1}, s_{2} \}$ of $E_{X}^{+}$.  Let $\ell = 3$, and consider the function $\alpha:E_{X}^{+} \rightarrow \mathbb{Z} \subseteq \mathbb{Z}_{3}$ given by
$$s_{1} \mapsto \alpha(s_{1}) = 1 \text{ and } s_{2} \mapsto \alpha(s_{2}) = 11. $$
By Corollary \ref{admissible corollary for bouquets} below, all the graphs $X(\mathbb{Z}/3^{n}\mathbb{Z},S,\alpha_{/n})$ are connected.  We obtain a $\mathbb{Z}_{3}$-tower above the bouquet graph $X$ which can be pictured as follows:
\begin{equation*}
\begin{tikzpicture}[baseline={([yshift=-1.7ex] current bounding box.center)}]
\node[draw=none,minimum size=2cm,regular polygon,regular polygon sides=1] (a) {};
\foreach \x in {1}
  \fill (a.corner \x) circle[radius=0.7pt];
\draw (a.corner 1) to [in=50,out=130,loop] (a.corner 1);
\draw (a.corner 1) to [in=50,out=130,distance = 0.5cm,loop] (a.corner 1);
\end{tikzpicture}
\longleftarrow \, \, \,
\begin{tikzpicture}[baseline={([yshift=-0.6ex] current bounding box.center)}]
\node[draw=none,minimum size=2cm,regular polygon,regular polygon sides=3] (a) {};

\foreach \x in {1,2,3}
  \fill (a.corner \x) circle[radius=0.7pt];

\path (a.corner 1) edge [bend left=15] (a.corner 2);
\path (a.corner 1) edge [bend right=15] (a.corner 2);
\path (a.corner 2) edge [bend left=15] (a.corner 3);
\path (a.corner 2) edge [bend right=15] (a.corner 3);
\path (a.corner 3) edge [bend left=15] (a.corner 1);
\path (a.corner 3) edge [bend right=15] (a.corner 1);
\end{tikzpicture}
\, \, \, \longleftarrow  \, \, \,
\begin{tikzpicture}[baseline={([yshift=-0.6ex] current bounding box.center)}]
\node[draw=none,minimum size=2cm,regular polygon,regular polygon sides=9] (a) {};

\foreach \x in {1,2,...,9}
  \fill (a.corner \x) circle[radius=0.7pt];
  
\foreach \y\z in {1/2,2/3,3/4,4/5,5/6,6/7,7/8,8/9,9/1}
  \path (a.corner \y) edge (a.corner \z);
  
\foreach \y\z in {1/3,2/4,3/5,4/6,5/7,6/8,7/9,8/1,9/2}
  \path (a.corner \y) edge  (a.corner \z);
\end{tikzpicture}
\, \, \, \longleftarrow \, \, \,
\begin{tikzpicture}[baseline={([yshift=-0.6ex] current bounding box.center)}]
\node[draw=none,minimum size=2cm,regular polygon,regular polygon sides=27] (a) {};

\foreach \x in {1,2,...,27}
  \fill (a.corner \x) circle[radius=0.7pt];
  
\foreach \y\z in {1/2,2/3,3/4,4/5,5/6,6/7,7/8,8/9,9/10,10/11,11/12,12/13,13/14,14/15,15/16,16/17,17/18,18/19,19/20,20/21,21/22,22/23,23/24,24/25,25/26,26/27,27/1}
  \path (a.corner \y) edge (a.corner \z);
  
\foreach \y\z in {1/12,2/13,3/14,4/15,5/16,6/17,7/18,8/19,9/20,10/21,11/22,12/23,13/24,14/25,15/26,16/27,17/1,18/2,19/3,20/4,21/5,22/6,23/7,24/8,25/9,26/10,27/11}
  \path (a.corner \y) edge (a.corner \z);
\end{tikzpicture}
\, \, \, \longleftarrow \ldots
\end{equation*}

\item Let $\ell=3$ again, and let now $X$ be the dumbbell multigraph
\begin{center}
\begin{tikzpicture}
\node [circle, draw, fill=white] (b1) at (1,1) {};
\node [circle, draw, fill=white] (b2) at (2,1) {};

\node (c) at (1.16,1) {};
\fill (c) circle[radius=0.7pt];

\node (d) at (1.84,1) {};
\fill (d) circle[radius=0.7pt];

\draw (b1) -- (b2);
\end{tikzpicture},
\end{center}
the section be
\begin{center}
\begin{tikzpicture}
\node [label = {\small $s_{1}$},circle, draw, fill=white] (b1) at (1,1) {};
\node [label = {\small $s_{3}$}, circle, draw, fill=white] (b2) at (2,1) {};

\node (c) at (1.16,1) {};
\fill (c) circle[radius=0.7pt];

\node (d) at (1.84,1) {};
\fill (d) circle[radius=0.7pt];

\draw[->=.6] (b1) -- (b2) node [midway,above] {\small $s_{2}$};
\end{tikzpicture},
\end{center}
and take the function $\alpha:S \rightarrow \mathbb{Z} \subseteq \mathbb{Z}_{3}$ given by
$$s_{1} \mapsto \alpha(s_{1}) = 1, s_{2} \mapsto \alpha(s_{2}) =0, \text{ and } s_{3} \mapsto \alpha(s_{3}) = 11. $$
By Theorem \ref{connectedness_cond} below, all the multigraphs $X(\mathbb{Z}/3^{n}\mathbb{Z},S,\alpha_{/n})$ are connected.  We obtain a $\mathbb{Z}_{3}$-tower above the dumbbell multigraph $X$ which can be pictured as follows:
\begin{equation*}
\begin{tikzpicture}[baseline={([yshift=-0.6ex] current bounding box.center)}]
\node [circle, draw, fill=white] (b1) at (1,1) {};
\node [circle, draw, fill=white] (b2) at (2,1) {};

\node (c) at (1.16,1) {};
\fill (c) circle[radius=0.7pt];

\node (d) at (1.84,1) {};
\fill (d) circle[radius=0.7pt];

\draw (b1) -- (b2);
\end{tikzpicture}
\, \, \, \longleftarrow \, \, \,
\begin{tikzpicture}[baseline={([yshift=-0.6ex] current bounding box.center)}]
\node[draw=none,minimum size=2cm,regular polygon,regular polygon sides=3] (a) {};
\node[draw=none, minimum size=1.3cm,regular polygon,regular polygon sides=3] (b) {};

\foreach \x in {1,2,...,3}
  \fill (a.corner \x) circle[radius=0.7pt];
  
\foreach \y in {1,2,...,3}
  \fill (b.corner \y) circle[radius=0.7pt];
  
\path (a.corner 1) edge (a.corner 2);
\path (a.corner 2) edge (a.corner 3);
\path (a.corner 3) edge (a.corner 1);

\path (a.corner 1) edge (b.corner 1);
\path (a.corner 2) edge (b.corner 2);
\path (a.corner 3) edge (b.corner 3);

\path (b.corner 1) edge (b.corner 3);
\path (b.corner 2) edge (b.corner 1);
\path (b.corner 3) edge (b.corner 2);

\end{tikzpicture}
\, \, \, \longleftarrow \, \, \,
\begin{tikzpicture}[baseline={([yshift=-0.6ex] current bounding box.center)}]
\node[draw=none,minimum size=2cm,regular polygon,regular polygon sides=9] (a) {};
\node[draw=none, minimum size=1.7cm,regular polygon,regular polygon sides=9] (b) {};

\foreach \x in {1,2,...,9}
  \fill (a.corner \x) circle[radius=0.7pt];
  
\foreach \y in {1,2,...,9}
  \fill (b.corner \y) circle[radius=0.7pt];
  
\foreach \x\z in {1/2,2/3,3/4,4/5,5/6,6/7,7/8,8/9,9/1}
  \path (a.corner \x) edge (a.corner \z);  
 
\foreach \x\z in {1/1,2/2,3/3,4/4,5/5,6/6,7/7,8/8,9/9} 
  \path (a.corner \x) edge (b.corner \z);
  
\foreach \x\z in {1/3,2/4,3/5,4/6,5/7,6/8,7/9,8/1,9/2} 
  \path (b.corner \x) edge (b.corner \z);
  
\end{tikzpicture}
\, \, \, \longleftarrow \, \, \,
\begin{tikzpicture}[baseline={([yshift=-0.6ex] current bounding box.center)}]
\node[draw=none,minimum size=2cm,regular polygon,regular polygon sides=27] (a) {};
\node[draw=none, minimum size=1.7cm,regular polygon,regular polygon sides=27] (b) {};

\foreach \x in {1,2,...,27}
  \fill (a.corner \x) circle[radius=0.7pt];
  
\foreach \y in {1,2,...,27}
  \fill (b.corner \y) circle[radius=0.7pt];
  
\foreach \x\z in {1/2,2/3,3/4,4/5,5/6,6/7,7/8,8/9,9/10,10/11,11/12,12/13,13/14,14/15,15/16,16/17,17/18,18/19,19/20,20/21,21/22,22/23,23/24,24/25,25/26,26/27,27/1}
  \path (a.corner \x) edge (a.corner \z);  
 
\foreach \x\z in {1/1,2/2,3/3,4/4,5/5,6/6,7/7,8/8,9/9,10/10,11/11,12/12,13/13,14/14,15/15,16/16,17/17,18/18,19/19,20/20,21/21,22/22,23/23,24/24,25/25,26/26,27/27} 
  \path (a.corner \x) edge (b.corner \z);
  
\foreach \x\z in {1/12,2/13,3/14,4/15,5/16,6/17,7/18,8/19,9/20,10/21,11/22,12/23,13/24,14/25,15/26,16/27,17/1,18/2,19/3,20/4,21/5,22/6,23/7,24/8,25/9,26/10,27/11}
  \path (b.corner \x) edge (b.corner \z);
  
\end{tikzpicture}
\, \, \, \longleftarrow \ldots
\end{equation*}
Note that the graphs $X(\mathbb{Z}/3^{n}\mathbb{Z},S,\alpha_{/n})$ are precisely the generalized Petersen graphs $G(3^{n},11)$.
\end{enumerate}

\par Let $u_X$ be the number of vertices of $X$ and label $V_X=\{v_1,\dots, v_{u_X}\}$. Let $o,t: E_X^+\rightarrow V_X$ be the \emph{origin and terminus functions} taking an edge to its source and target vertex, respectively. Thus, we have that $i(e)=(o(e), t(e))$. For $v\in V_X$, set $E_v:=\{e\in E_X^+\mid o(e)=v\}$, and let $\op{val}_X(v):=\# E_v$ be the \emph{valency} of $v$. The divisor group $\op{Div}(X)$ is the free abelian group on the vertices $V_X$ of $X$. It consists of formal sums of the form $D=\sum_{v\in V_X} n_v v$, where $n_v\in \Z$ for all $v$. The degree of $D$ is the sum $\op{deg}(D):=\sum_v n_v\in \Z$, and this defines a homomorphism 
\[\op{deg}:\op{Div}(X)\rightarrow \Z,\] the kernel of which is denoted $\op{Div}^0(X)$. Let $\mathcal{M}(X)$ be the abelian group of $\Z$-valued functions on $V_X$.  
 The discrete Laplacian map ${\rm div}:\mathcal{M}(X) \to {\rm Div}(X)$ is defined as follows.  For $f \in \mathcal{M}(X)$, one sets
$${\rm div}(f) = -\sum_{v}m_{v}(f) \cdot v,$$
where
$${m_{v}(f)} = \sum_{e \in E_{X,v}}\left(f(t(e)) - f(o(e))\right).$$ We let $\op{Pr}(X)$ be the image of $\op{div}$. It is easy to see that $\op{Pr}(X)$ is a subgroup of $\op{Div}^0(X)$, since the subjacent directed graph is eulerian, meaning that every vertex has in-degree equals to its out-degree. The Picard group is taken to be the quotient $\op{Pic}(X):=\op{Div}(X)/\op{Pr}(X)$, and we set $\op{Pic}^0(X):=\op{Div}^0(X)/\op{Pr}(X)$. Set $\kappa_X$ to be the cardinality of $\op{Pic}^0(X)$, which is a finite group and is equal to the number of spanning trees in $X$.

\par The \emph{valency matrix} is the $u_X\times u_X$ diagonal matrix $D=(d_{i,j})$, with 
\[d_{i,j}=\begin{cases}
\op{val}_X(v_i) & \text{ for }i=j,\\
0 & \text{ for }i\neq j.
\end{cases}\]
The \emph{adjacency matrix} $A=(a_{i,j})$ is given by
\[a_{i,j}=\begin{cases}
\text{Twice the number of undirected loops at }i,\text{  } & i=j;\\
\text{The number of undirected edges connecting the }i\text{-th and }j\text{-th vertices}, & i\neq j.
\end{cases}\]
The matrix $Q:=D-A$ is called the \emph{Laplacian matrix}.

\par We let $\Z_\ell[x]$ (resp. $\Z_\ell\llbracket T \rrbracket $) denote the polynomial ring (resp. formal power series ring) in the indeterminate $x$ (resp. $T$) with coefficients in $\Z_\ell$. Denote by $\Z_\ell[x;\Z_\ell]$ the ring whose elements are expressions of the form $f(x)=\sum_{a} c_a x^a$, where $a$ runs over a finite set in $\Z_\ell$ and the coefficient $c_a$ is in $\Z_\ell$. Here, multiplication is defined by setting $x^ax^b=x^{a+b}$ for $a,b\in \Z_\ell$. Associate a matrix $M(x)$ to the pair $(X,\alpha)$ as follows
\begin{equation}\label{M matrix}M(x)=M_{X,\alpha}(x):=D - \left(\sum_{\substack{s \in S \\ {i}(s) = (v_{i},v_{j})}} x^{\alpha(s)} +  \sum_{\substack{s \in S \\ { i}(s) = (v_{j},v_{i})}} x^{-\alpha(s)} \right) \in M_{u \times u}(\mathbb{Z}[x;\mathbb{Z}_{\ell}]).\end{equation}For $b\in \Z_\ell$, let ${b\choose n}\in \Z_\ell$ be given by
\[{b\choose n}:=\frac{b(b-1)\dots (b-n+1)}{n!}.\]Let $(1+T)^b$ be the formal power series 
\[(1+T)^b:=\sum_{n=0}^\infty {b\choose n} T^n.\] With this convention in place, we define the characteristic series to be \begin{equation}\label{char series def} f(T)=f_{X,\alpha}(T)=\op{det}M(1+T),\end{equation} which is a formal power series in $\Z_\ell\llbracket T\rrbracket $. The convention used in this article differs from that used in \cite{mcgownvallieresIII}, where $-T$ is used in place of $T$. This convention simplifies our calculations.

\begin{lemma}\label{f(0)=0 lemma}
Let $X$ be a multigraph, $\ell$ be a prime and $\alpha:S\rightarrow \Z_\ell$ a function. The constant term of the characteristic series $f(T)$ (see \eqref{char series def}) is equal to $0$.
\end{lemma}
\begin{proof}
Note that the constant term of $f(T)$ is $\op{det}(Q)$. The vector $(1,\ldots,1)^t$ is in the null-space of $Q$, see \cite[proof of Corollary~9.10]{CP-book} for further details. Therefore, $Q$ is singular and the result follows.
\end{proof}


The following result shows that $f_{X, \alpha}(T)\neq 0$ under the assumptions imposed on $\alpha$. The proof is essentially contained in \cite{Vallieres:2021, mcgownvallieresII, mcgownvallieresIII}, we recall these details. First, we introduce some further notation. Let $Y/X$ be an abelian cover of multigraphs and set $G:=\op{Gal}(Y/X)$. Let $\widehat{G}$ be the group of characters $\psi: G\rightarrow \mathbb{C}^\times$. Given $\psi\in \widehat{G}$, denote the Artin--Ihara L-function by
\[L_X(z, \psi):=\prod_{\mathfrak{c}}\left(1-\psi\left(\left(\frac{Y/X}{\frak{c}}\right)\right)z^{l(\mathfrak{c})}\right)^{-1}.\]Here, the product runs over all primes $\mathfrak{c}$ of $X$, $l(\mathfrak{c})$ is the length of $\mathfrak{c}$ and $\left(\frac{Y/X}{\frak{c}}\right)$ is the Frobenius automorphism at $\mathfrak{c}$. (Recall that a prime in a multigraph is an equivalence class of closed paths that have no backtrack or tail, and that are not multiples of a strictly smaller closed path, where two such paths are identified if they are the same up to their starting vertex. For further details, the reader is referred to \cite[Chapter~18]{Terras:2011}.) Let $\chi(X)$ denote the Euler characteristic of $X$, which is defined as follows
\[\chi(X):=|V_X|-|E_X|.\] We shall impose the following assumption without further mention.

\begin{assumption}\label{euler char assumption}
Throughout, we assume that $\chi(X)\neq 0$. 
\end{assumption}

That $\chi(X)\neq 0$ means that $X$ does not consist of a single cycle, or a "cycle with hair", cf. \cite[Figure 2.1]{Terras:2011}. When $\chi(X)=0$, $X$ is referred to as a "bad graph", cf. P.10 of \emph{op. cit.} (and references therein), since many arithmetic properties of zeta functions are not satisfied in this case. For instance, the analog of the prime number theorem for such bad graphs does not hold, cf. Theorem 10.1 of \emph{op. cit.}

\begin{definition} \label{artinian}
Let $Y/X$ be an abelian cover of multigraphs with automorphism group $G$, and for each $i=1,\ldots,u_{X}$, let $w_{i}$ be a fixed vertex in the fiber of $v_{i}$.
\begin{enumerate}
\item For $\sigma \in G$, we let the matrix $A(\sigma)$ to be the $u_{X} \times u_{X}$ matrix $A(\sigma)=(a_{i,j}(\sigma))$ defined via
\begin{equation*}
a_{i,j}(\sigma) = 
\begin{cases}
\text{Twice the number of undirected loops at the vertex }w_{i}, &\text{ if } i=j \text{ and } \sigma = 1;\\
\text{The number of undirected edges connecting $w_{i}$ to $w_{j}^{\sigma}$}, &\text{ otherwise}.
\end{cases}
\end{equation*}
\item For $\psi \in \widehat{G}$, we set
$$A_{\psi} = \sum_{\sigma \in G} \psi(\sigma) \cdot A(\sigma). $$
\end{enumerate}
\end{definition}

\begin{lemma}\label{f(T) is not zero}
Let $\alpha:S\to \mathbb{Z}_{\ell}$ be a function for which Assumption \ref{main assumption} is satisfied, then $f_{X,\alpha}(T) \neq 0$. 
\end{lemma}
\begin{proof}
The relationship between the power series $Q(T)$ of \cite[Theorem 6.1]{mcgownvallieresIII} and $f_{X,\alpha}(T)$ is $Q(T) = f_{X,\alpha}(-T)$.  If $Y/X$ is an abelian cover of multigraphs and $\psi$ is a character of $G = {\rm Gal}(Y/X)$, then the three-term determinant formula \cite[Theorem 18.15]{Terras:2011} for the Artin-Ihara $L$-function gives
$$L_{X}(z,\psi)^{-1} = (1-z^{2})^{-\chi(X)} \cdot {\rm det}(I - A_{\psi}z + (D-I)z^{2}),$$
where $A_{\psi}$ is the twisted adjacency matrix of Definition \ref{artinian}, $D$ is the valency matrix of $X$, and $\chi(X)$ is the Euler characteristic of $X$. From now on, we let
$$h_{X}(z,\psi) =  {\rm det}(I - A_{\psi}z + (D-I)z^{2}) \in \mathbb{Z}[\psi][z].$$
Coming back to the abelian $\ell$-tower
$$X  \leftarrow X(\mathbb{Z}/\ell\mathbb{Z},S,\alpha_{/1}) \leftarrow X(\mathbb{Z}/\ell^{2}\mathbb{Z},S,\alpha_{/2}) \leftarrow \ldots \leftarrow X(\mathbb{Z}/\ell^{n}\mathbb{Z},S,\alpha_{/n}) \leftarrow \ldots$$
associated to $\alpha$, for all positive integer $n$, \cite[Corollary 5.6]{mcgownvallieresIII} implies
\begin{equation} \label{spe_pow1}
Q(1- \zeta_{\ell^{n}}) = h_{X}(1,\psi_{n}), 
\end{equation}
where $\psi_{n}$ is the character of $\mathbb{Z}/\ell^{n}\mathbb{Z} \simeq {\rm Gal}(X(\mathbb{Z}/\ell^{n}\mathbb{Z},S,\alpha_{/n})/X)$ satisfying $\psi_{n}(\bar{1}) = \zeta_{\ell^{n}}$.  Here, it is understood that $\zeta_{\ell^{n}} = \exp(2 \pi i/\ell^{n})$ and that an embedding of $\overline{\mathbb{Q}}$ into $\overline{\mathbb{Q}}_{\ell}$ has been fixed once and for all so that we view $\zeta_{\ell^{n}}$ and the values of the characters $\psi_{n}$ inside $\overline{\mathbb{Q}}_{\ell}$.  From \cite[\S 3]{Vallieres:2021}, one has
$$\ell^{n} \cdot \kappa_{n} = \kappa_{X} \cdot \prod_{\psi \neq \psi_{0}} h_{X}(1,\psi), $$
where the product is over all non-trivial characters of $\mathbb{Z}/\ell^{n}\mathbb{Z}$, and $\kappa_n$ denotes the number of spanning trees of $X(\mathbb{Z}/\ell^{n}\mathbb{Z},S,\alpha_{/n})$.  
The above formula crucially requires the Assumption \ref{euler char assumption}. It follows that
$$h_{X}(1,\psi) \neq 0 $$
for all non-trivial characters, and thus $Q(T) \neq 0$ by (\ref{spe_pow1}).
\end{proof}

\par Associated to the characteristic series, are the Iwasawa $\mu$ and $\lambda$ invariants, which we now proceed to define. A polynomial is said to be a \emph{distinguished polynomial} if it is monic and all non-leading coefficients are divisible by $\ell$.
\begin{definition}\label{defn of mu and lambda}
By the Weierstrass Preparation theorem, we may factor the characteristic polynomial as follows $f(T)=\ell^{\mu} g(T) u(T)$, where $g(T)$ is a distinguished polynomial in $\Z_\ell\llbracket T\rrbracket $, and $u(T)$ is a unit in $\Z_\ell\llbracket T\rrbracket $. Since $u(T)$ is a unit in $\Z_\ell\llbracket T\rrbracket$, the constant term $u(0)$ is a unit in $\Z_\ell$. According to Lemma \ref{f(0)=0 lemma}, the constant term of $f(T)$ is equal to $0$, and hence $T$ divides $g(T)$. In particular, we have that $\op{deg} g\geq 1$. The number $\mu\in \Z_{\geq 0}$ is the $\mu$-invariant and depends on $(\ell,\alpha)$. Throughout it is assumed that Assumption \ref{main assumption} is satisfied for $\alpha$. Thus, we denote this invariant by $\mu_{\ell}(X,\alpha)$, however, when $\ell$, $X$ and $\alpha:S\rightarrow \Z_\ell$ are understood, we simply denote it by $\mu$. On the other hand, the $\lambda$-invariant is defined to be the quantity 
\[\lambda=\lambda_{\ell}(X,\alpha):=\op{deg} g(T)-1.\] 
\end{definition}

\begin{remark}
\hfill
\begin{itemize}
    \item According to Lemma \ref{f(T) is not zero}, we have that $f(T)\neq 0$, and thus the $\mu$ and $\lambda$ invariants are defined. 
    \item The $\nu$-invariant is not detectable from the characteristic series, and we postpone its definition till Definition \ref{def of nu}.
    \item In our setting the $\lambda$-invariant is \emph{not} the degree of $g(T)$, i.e., the conventional definition of the $\lambda$-invariant of a power series $f(T)$. The shift by $-1$ is necessary in studying the growth of Picard numbers in a $\Z_\ell$-tower.
\end{itemize}

\end{remark}
\par Let $\alpha:S\rightarrow \Z_\ell$ be a function 
such that Assumption \ref{main assumption} is satisfied. Given $n\geq 0$, let $X_n=X(\Z/\ell^n\mathbb{Z}, S, \alpha_{/n})$ and write  $\kappa_n$ for the number of spanning trees of $X_n$. Let $|\cdot|_\ell:\Q_\ell\rightarrow \mathbb{R}_{\geq 0}$ 
be the absolute value normalized by setting $|\ell|_\ell^{-1}=\ell$. The $\ell$-part of $\kappa_n$ is defined as $\kappa_\ell(X_n):=|\kappa_n|_\ell^{-1}$. The following result is \cite[Theorem 6.1]{mcgownvallieresIII}.
\begin{theorem}[McGown-Valli\`{e}res]
Let $\ell$ be a prime and $\alpha$ be as above. Then, there exists $n_0>0$ and $\nu\in \Z$, such that
\begin{equation}\label{asymptotic formula}\kappa_\ell(X_n)=\ell^{(\mu \ell^n+\lambda n+\nu)},\end{equation} for all $n\geq n_0$. The invariants $\mu$ and $\lambda$ are given as in Definition \ref{defn of mu and lambda}.
\end{theorem}

\begin{definition}\label{def of nu}
The $\nu$-invariant is the uniquely determined constant for which the asymptotic formula \eqref{asymptotic formula} holds.
\end{definition}

\subsection{Connectedness of multigraphs in abelian $\ell$-towers}\label{connectedness assumption section}
Let $X$ be a connected multigraph, $\ell$ be a prime number and let $\alpha:S\rightarrow \Z_\ell$. In this section, we introduce explicit conditions on $X$ which are sufficient for Assumption \ref{main assumption} to hold. Let
$$X=X_0\leftarrow X_1\leftarrow X_2\leftarrow \dots \leftarrow X_n \leftarrow \dots $$
be an abelian $\ell$-tower over $X$.  Set also $G_{n} = {\rm Gal}(X_{n}/X)$.  We have natural group morphisms
$$\gamma_{n}:H_{1}(X,\mathbb{Z}) \to G_{n}$$
that are compatible so that we get a group morphism
$$\gamma:H_{1}(X,\mathbb{Z}) \to \varprojlim_n G_{n} \simeq \mathbb{Z}_{\ell}.$$
Fix now a spanning tree of $X$, say $\mathfrak{T}\subseteq E_X$,  
a section $\gamma:E_{X} \to E_{X}^+$, and set $S:=\gamma(E_{X})$.  We let $S(\mathfrak{T})$ be the collection of $s \in S$ such that $\pi(s) \in \mathfrak{T}$.  For each $s \in S$, let $c_{s} = d_{s} \cdot s$ be the closed path in $X$, where $d_{s}$ is the unique geodesic path in $\mathfrak{T}$ going from $t(s)$ to $o(s)$.  Let $\langle c_{s} \rangle$ be the corresponding cycle in $H_{1}(X,\mathbb{Z})$.  Note that if $s \in S(\mathfrak{T})$, then $\langle c_{s} \rangle = 0$.  
The collection
$$\{ \langle c_{s} \rangle \, | \, s \in S \smallsetminus S(\mathfrak{T}) \}$$
forms a basis of $H_{1}(X,\mathbb{Z})$ (cf. \cite[\S 4.3]{Sunada:2013}), so that in particular
$${\rm rank}_{\mathbb{Z}}(H_{1}(X,\mathbb{Z})) = |S \smallsetminus S(\mathfrak{T})| = |E_{X}| - |V_{X}| + 1.$$
Define now $\alpha:S \to \mathbb{Z}_{\ell}$ via $s \mapsto \alpha(s) = \gamma(\langle c_{s} \rangle)$. By successively lifting the tree $\mathfrak{T}$ to $X_{1}$, then to $X_{2}$ and so on, one can construct isomorphisms of multigraphs
$$\phi_{n}:X(\mathbb{Z}/\ell^{n},S,\alpha_{/n}) \xrightarrow{\sim} X_{n}$$
for which all the squares and the triangle in the diagram
\begin{equation*} \label{nice_dia}
\begin{tikzcd}[sep=1.8em, font=\small]
         & \arrow[ld] X(\mathbb{Z}/\ell\mathbb{Z},S,\alpha_{/1})  \arrow[dd, "\phi_{1}"]  & \arrow[l]X(\mathbb{Z}/\ell^{2}\mathbb{Z},S,\alpha_{/2})   \arrow[dd, "\phi_{2}"] &\arrow[l] \ldots  &\arrow[l] \arrow[dd, "\phi_{n}"]  X(\mathbb{Z}/\ell^{n}\mathbb{Z},S,\alpha_{/n}) & \arrow[l]  \ldots\\  
         X &  & & & & \\
          & \arrow[lu]X_{1}  & \arrow[l]X_{2}   & \arrow[l] \ldots  & \arrow[l] X_{n}  \arrow[l] & \arrow[l] \ldots
\end{tikzcd}
\end{equation*}
commute.  In other words, once a spanning tree $\mathfrak{T}$ for $X$ has been fixed (which also specifies a $\mathbb{Z}$-basis for $H_{1}(X,\mathbb{Z})$), every abelian $\ell$-tower over $X$ arises from a function $\alpha:S \to \mathbb{Z}_{\ell}$ satisfying $\alpha(s) = 0$ for all $s \in S(\mathfrak{T})$.

\begin{theorem} \label{connectedness_cond}
Let $X$ be a connected multigraph and let $\mathfrak{T}\subseteq E_X$ be a fixed spanning tree of $X$.  Furthermore, let $S$ be the image of a section as above, and let $\alpha:S \to \mathbb{Z}_{\ell}$ be a function satisfying $\alpha(s) = 0$ for all $s \in S(\mathfrak{T})$.  Then $X(\mathbb{Z}/\ell^{n}\mathbb{Z},S,\alpha_{/n})$ is connected for all positive integer $n$ if and only if there exists $s \in S \smallsetminus S(\mathfrak{T})$ such that $\alpha(s)\in \Z_\ell^\times$.
\end{theorem}
\begin{proof}
Let us assume first that there exists $s \in S \smallsetminus S(\mathfrak{T})$ such that $(\alpha(s),\ell)=1$.  Let also $X_{n} = X(\mathbb{Z}/\ell^{n}\mathbb{Z},S,\alpha_{/n})$. Let $v\in V_X$, $e\in E_X^+$ and $\bar{a}\in \Z/\ell^n\Z$. The map $p_{n}:X_{n} \to X$ defined via
$$p_{n}\left((v,\bar{a})\right) = v \text{ and } p_{n}\left((e,\bar{a})\right) = e  $$
is a morphism of multigraphs (recall that the vertices and (directed) edges of $X_n$ consist of elements in $V_X\times \Z/\ell^n\Z$ and $E_X^+\times \Z/\ell^n\Z$, respectively). Furthermore, for all $w \in V_{X_{n}}$, the restriction $p_{n}|_{E_{X_{n},w}}$ induces a bijection
$$p_{n}|_{E_{X_{n}},w}: E_{X_{n},w} \xrightarrow{\sim} E_{X,p_{n}(w)}.$$
It follows that every path $c$ in $X$ has a unique lift $\tilde{c}$ to $X_{n}$ once the starting point in the fiber of $o(c)$ is fixed.  Furthermore, since $p_{n}$ is a morphism of multigraphs one has
\begin{equation} \label{lift}
p_{n}(o(\tilde{c})) = o(c) \text{ and } p_{n}(t(\tilde{c})) = t(c). 
\end{equation}
Now, $\mathbb{Z}/\ell^{n}\mathbb{Z}$ acts on $X_{n}$ via
$$\bar{a} \cdot (v,\bar{b}) = (v, \bar{a} + \bar{b}) \text{ and } \bar{a} \cdot (e,\bar{b}) = (e,\bar{a} + \bar{b}), $$
which gives an injective group morphism
$$\mathbb{Z}/\ell^{n}\mathbb{Z} \hookrightarrow {\rm Aut}_{p_{n}}(X_{n}/X). $$
Let $v_{1},v_{2} \in V_{X}$ be arbitrary.  Since $X$ is connected there is a path $c$ in $X$ satisfying $o(c)=v_{1}$ and $t(c) = v_{2}$.  Thus, by lifting $c$ and using (\ref{lift}) above, for any vertex $w \in p_{n}^{-1}(v_{1})$, there exists a path $\tilde{c}$ in $X_{n}$ going from $w$ to another vertex in $p_{n}^{-1}(v_{2})$.  It follows that in order to show that $X_{n}$ is connected, it suffices to show that every vertex in a single fiber can be joined via a path in $X_{n}$.  Under our assumption, there exists $s_{0} \in S \smallsetminus S(\mathfrak{T})$ such that $(\alpha(s_{0}),\ell)=1$, and thus $\alpha_{/n}(s_{0})$ is a generator for $\mathbb{Z}/\ell^{n}\mathbb{Z}$.  Let $v \in V_{X}$ be such that $v = t(s_{0})$ and let $\sigma_{0} = \alpha_{/n}(s_{0}) \in \mathbb{Z}/\ell^{n}\mathbb{Z}$.  The closed path $c_{s_{0}} = d_{s_{0}}\cdot s_{0}$ can be lifted to a path $\tilde{c}_{s_{0}}$ in $X_{n}$ starting at $(v,\bar{0})$.  Furthermore, since $\alpha(s) = 0$ for all $s \in S(\mathfrak{T})$, we have $t(\tilde{c}_{s_{0}}) = (v,\sigma_{0})$.  The path $\tilde{c}_{s_{0}} \cdot \sigma_{0}(\tilde{c}_{s_{0}})$ will start at $(v,\bar{0})$ and end at $(v,2\cdot \sigma_{0})$.  Keeping going like this, we can find a path in $X_{n}$ going from $(v,\bar{0})$ to $(v,m\cdot \sigma_{0})$ for all positive integer $m$.  Since $\sigma_{0}$ generates $\mathbb{Z}/\ell^{n}\mathbb{Z}$, this shows that every vertex in the fiber of $v$ can be joined via a path in $X_{n}$.  It follows that $X_{n}$ is connected. 

Conversely, assume that all multigraphs $X_{n}$ are connected.  In particular, $X_{1}$ is connected.  To simplify the notation, write $S \smallsetminus S(\mathfrak{T}) = \{s_{1},\ldots,s_{t} \}$.  The function $\alpha$ induces a function on $C_{1}(X,\mathbb{Z})$ and thus also on $H_{1}(X,\mathbb{Z})$ by restriction.  Let also $v \in V_{X}$ be a fixed vertex.  Since we are assuming that $X_{1}$ is connected, there exists a path $\tilde{c}$ in $X_{1}$ going from $(v,\bar{0})$ to $(v,\bar{1})$.  The path $c = p_{1}(\tilde{c})$ is then a closed path in $X$, and $\alpha(\langle c \rangle) = \bar{1}$.  Since the $\langle c_{s_{i}} \rangle$ for $i=1,\ldots,t$ form a $\mathbb{Z}$-basis of $H_{1}(X,\mathbb{Z})$, we have
$$\langle c \rangle = n_{1} \langle c_{s_{1}} \rangle + \ldots + n_{t} \langle c_{s_{t}}\rangle, $$
for some $n_{i} \in \mathbb{Z}$.  Applying $\alpha_{/1}$ to this last equation gives
$$\bar{1} = n_{1} \alpha_{/1}(s_{1}) + \ldots + n_{t}\alpha_{/1}(s_{t}). $$
Now, if all the $\alpha(s_{i})$ were divisible by $\ell$, then we would get a contradiction. 
\end{proof}
The following corollary is most likely well known, we list it here for convenience.
\begin{corollary}\label{admissible corollary for bouquets}
Let $X$ be a bouquet and $\alpha:S\rightarrow \Z_\ell$. Then, Assumption \ref{main assumption} is satisfied if and only if there exists $s\in S$ for which $\ell \nmid \alpha(s)$.  
\end{corollary}
\begin{proof}
Since $X$ is a bouquet, the only possible spanning tree consists of the single vertex of $X$ (and no edges). The result follows from Theorem \ref{connectedness_cond}.
\end{proof}
\subsection{Basic properties of Iwasawa invariants}
\par In this subsection, we prove a number of basic results about the Iwasawa $\mu$ and $\lambda$ invariants associated to multigraphs. The graph $X$ is called a \emph{bouquet} if there is $1$ vertex, in which case $M(x)$ is a $1\times 1$-matrix. Suppose that there are $t$ loops 
in total, labelled $e_1,\dots, e_t$ and $\alpha=(\alpha_1,\dots, \alpha_t)$ sends $e_i\mapsto \alpha_i$. We denote the corresponding bouquet by $X_t$. We shall set $\mu:=\mu_\ell(X_t;\alpha)$ and $\lambda:=\lambda_\ell(X_t;\alpha)$. Given $a\in \Z_\ell$, let \[f_a(T)=2-(1+T)^a-(1+T)^{-a}=-\sum_{i\geq 2}\left({a\choose i}+{-a\choose i}\right) T^i.\] We have that $f(T)=\sum_{j=1}^t f_{\alpha_j}(T)$, and therefore we may write
\[f(T)=\sum_{n\geq 2} \beta_{n} T^n,\] with $\beta_2=-\left(\sum_{j} \alpha_j^2\right)$. In particular, we find that $\op{ord}_{T=0} f(T)\geq 2$, and thus $\lambda\geq 1$. It is easy to see that $\mu=0$ and $\lambda=1$ if and only if $\beta_2$ is a unit in $\Z_\ell$. These are the minimal values of the Iwasawa invariants attainable by a bouquet. The following result shows that the Iwasawa invariants $\mu$ and $\lambda$ can be arbitrarily large.
\begin{lemma} \label{arb_large}
Let $n_1\in \Z_{\geq 0}$, $n_2\in \Z_{\geq 1}$, and let $\ell$ be a prime number. There is a bouquet $X_t$ consisting of $t = \ell^{n_{1}+1} + \ell^{n_{1}}$ loops and a choice of vector $\alpha\in \left(\Z_{\ell}^t\backslash (\ell \Z_{\ell})^t\right)$ (which depends on both $n_1$ and $n_2$) such that the associated Iwasawa invariants satisfy
\[\mu_{\ell}(X_t, \alpha)= n_1, \quad\lambda_{\ell}(X_t, \alpha)= 2\ell^{n_2}-1.\]
\end{lemma}
\begin{proof}
Set $t=\ell^{n_1+1}+\ell^{n_1}$ and let $\alpha=(\alpha_1, \dots, \alpha_t)$ with 
\[\alpha_i=\begin{cases}
 1&\text{ if }i\leq \ell^{n_1+1},\\
 \ell^{n_2} &\text{ if }i> \ell^{n_1+1}.\\
\end{cases}\]
We have that \[f_1(T)=2-(1+T)-(1+T)^{-1}=-(1+T)^{-1}T^2,\] and that 
\[f_{\ell^{n_2}}(T)=2-(1+T)^{\ell^{n_2}}-(1+T)^{-\ell^{n_2}}=-(1+T)^{-\ell^{n_2}}\left(1-(1+T)^{\ell^{n_2}}\right)^2.\]As a result, we find that \[\begin{split}f(T)=& \ell^{n_1}\left(\ell f_1(T)+f_{\ell ^{n_2}}(T)\right),\\
= & -\ell^{n_1} (1+T)^{-\ell ^{n_2}}\left(\ell (1+T)^{\ell ^{n_2-1}}T^2 +\left(1-(1+T)^{\ell ^{n_2}}\right)^2\right).\end{split}\] The degree of $\left(1-(1+T)^{\ell ^{n_2}}\right)^2$ is $2 \ell ^{n_2}$, which is larger than the degree of $(1+T)^{\ell ^{n_2-1}}T^2$. Furthermore, all non-leading coefficients of $\left(1-(1+T)^{\ell ^{n_2}}\right)^2$ are divisible by $\ell$. We have shown that $P(T):=\ell (1+T)^{\ell ^{n_2-1}}T^2 +\left(1-(1+T)^{\ell ^{n_2}}\right)^2$ is a distinguished polynomial. Thus, $f(T)=\ell^{n_1}u(T)P(T)$, where $u(T)=-(1+T)^{-\ell^{n_2}}$ is a unit in the Iwasawa algebra and $P(T)$ is a distinguished polynomial of degree $2\ell^{n_2}$. Thus, we have shown that $\mu=n_1$ and $\lambda=2\ell ^{n_2}-1$. 
\end{proof}

\begin{theorem}\label{lambda is odd}
Let $X$ be a multigraph, let $\ell$ be an odd prime and $\alpha:S\rightarrow \Z_\ell$ be a function such that Assumption \ref{main assumption} is satisfied. Then, $\lambda_\ell(X, \alpha)$ is odd. 
\end{theorem}

\begin{proof}
Set $\dot{T}:=(1+T)^{-1}-1$ and note that \begin{equation}\label{above equation T dot}(1+\dot{T})^{a}+(1+\dot{T})^{-a}=(1+T)^{-a}+(1+T)^a\end{equation}for all $a\in\Zl$. Note that $M(x)^t=M(x^{-1})$ and
it follows from \eqref{above equation T dot} that
\[f(\dot{T})=\op{det}M(1+\dot{T})=\op{det}M\left((1+T)^{-1}\right)=\op{det}M(1+T)=f(T).\]

Write $f(T)=\ell^\mu P(T) u(T)$, where $P(T)=\sum_{j=0}^{\lambda} \xi_j T^j+ T^{\lambda+1}$ and $\xi_j$ are all divisible by $\ell$. Note that \[(1+T)^{\lambda+1} P(\dot{T})=\sum_{j=0}^{\lambda} (-1)^{j}\xi_j (1+T)^{\lambda+1-j}T^j+ (-1)^{\lambda+1} T^{\lambda+1}.\] 
We see that $P(-1)=\sum_{j=0}^{\lambda} \xi_j (-1)^j+(-1)^{\lambda+1}\equiv (-1)^{\lambda+1}\mod{\ell}$, and hence, $P(-1)$ is a unit in $\Z_\ell$. The above relation shows that $P(-1)^{-1}(1+T)^{\lambda+1}P(\dot{T})$ is the distinguished polynomial dividing $f(\dot{T})=f(T)$, thus
\[P(T)=P(-1)^{-1}(1+T)^{\lambda+1}P(\dot{T})\] and 
\[u(T)=P(-1)(1+T)^{-\lambda-1} u(\dot{T}).\]
The latter relation specialized to $T=0$ gives $P(-1)=1$. Thus, we have shown that $1\equiv(-1)^{\lambda+1}\mod{\ell}$, and since $\ell$ is assumed to be odd, it follows that $\lambda$ is odd.
\end{proof}

\section{Equations over finite fields}

\par Let $\ell$ be a prime ($\ell=2$ is allowed) and $q$ be a power of $\ell$. Denote by $\F_q$ the finite field with $q$ elements, and let $f_1, \dots, f_m\in \F_q[x_1, \dots, x_n]$. In this section, we catalogue some well known results on the number of solutions to the system of equations $f_1=f_2=\dots =f_m=0$.
\subsection{Systems of equations}
\par We begin with a generalization of the Chevalley--Warning theorem to a system of equations in many variables. Associate to $I=(i_1, \dots, i_n)\in \Z_{\geq 0}^n$ the monomial $x_I:=x_1^{i_1}x_2^{i_2}\dots x_n^{i_n}$, then define the degree to be the sum \[\op{deg}\left(x_I\right):=i_1+i_2+\dots +i_n.\]
Given $f(x_1, \dots, x_n)=\sum_I c_I x_I$, the degree of $f$ is defined as follows \[\op{deg} f:=\op{max}\left\{\op{deg}(x_I)\mid c_I\neq 0\right\}.\] The polynomial is homogenous of degree $d$ if $\op{deg}(x_I)=d$ for all $I$ such that $c_I\neq 0$. The following is often referred to as \emph{Warning's second theorem}, cf. \cite{warning1936bemerkung}.
\begin{theorem}\label{lower bound on number of solns}
Let $m, n\in \Z_{\geq 1}$ and $f_1,\dots, f_m\in \F_q[x_1, \dots, x_n]$. Assume that 
\[d:=\sum_{i=1}^m \op{deg}(f_i)<n\] and that for all $i$, we have that $f_i(0,0,\dots, 0)=0$. Let $N$ be the number of solutions of $f_1=f_2=\dots=f_m=0$. Then, $N$ is divisible by $\ell$ and $N\geq q^{n-d}$.
\end{theorem}
\begin{proof}
The result follows from \cite[Theorem 6.11]{lidl1997finite}.
\end{proof}


\subsection{Quadratic forms and their solutions}
\par A quadratic form is a homogenous polynomial of degree $2$. Assume throughout that $q$ is odd, let $\eta:\F_q^\times\rightarrow \{\pm 1\}$ be the quadratic character of $\F_q$, and extend $\eta$ to a function on $\F_q$ by setting $\eta(0)=0$. Given a quadratic form $Q(x_1,\dots, x_n)$ over $\F_q$, express $Q$ as 
\[Q(x_1,\dots, x_n)=\sum_{i,j} a_{i,j} x_i x_j,\] with $a_{i,j}=a_{j,i}$. The rank of $Q$ is the rank of the associated matrix $A=(a_{i,j})$ over $\F_q$. We let $\op{det} Q$ denote the determinant of $A$. The form $Q$ is non-degenerate if $\op{rank}A=n$, i.e., $A$ is non-singular. On the other hand, $Q$ is diagonal if $a_{i,j}=0$ unless $i=j$.

\begin{proposition}\label{diagonal form presentation}
Suppose that $\ell$ is odd and $q$ is a power of $\ell$. Then, any quadratic form over $\F_q$ is equivalent to a diagonal form.
\end{proposition}
\begin{proof}
The above result is \cite[Theorem 6.21]{lidl1997finite}.
\end{proof}

Let $v:\F_q\rightarrow \Z$ be the function defined by 
\[v(b):=\begin{cases}
-1, &\text{ if }b\neq 0; \\
(q-1),&\text{ if }b=0.
\end{cases}\]

\begin{theorem}\label{thm:quadratic-Fq}
Let $Q(x_1, \dots, x_n)$ be a non-zero quadratic form over $\F_q$, and assume that $q$ is odd and $n$ is even. For $b\in \F_q$, let $N\left(Q(x_1,\dots, x_n)=b\right)$ be the number of solutions to $Q(x_1,\dots, x_n)=b$. Then, we have that 
\[N\left(Q(x_1,\dots, x_n)=b\right)= q^{n-1}+ v(b)q^{\frac{n-2}{2}}\eta\left((-1)^{n/2} \Delta\right)\]
where $\eta$ is the quadratic character of $\F_q$ and $\Delta=\op{det} Q$. 
\end{theorem}
\begin{proof}
The result above is \cite[Theorem 6.26]{lidl1997finite}.
\end{proof}
Given a degenerate quadratic form $Q(x_1, \dots, x_n)$, we let $\mathfrak{r}$ denote the rank of $Q$. According to Proposition \ref{diagonal form presentation}, $Q$ is equivalent to the form $\sum_{i=1}^{\mathfrak{r}} a_i x_i^2$, where $a_i$ are all non-zero. We set $Q_0(x_1,\dots, x_{\mathfrak{r}})\in \F_q[x_1, \dots, x_{\mathfrak{r}}]$ to be the nondegenerate form $Q_0(x_1, \dots, x_{\mathfrak{r}})=\sum_{i=1}^{\mathfrak{r}} a_i x_i^2$, with discriminant $\Delta(Q_0)=\prod_{i=1}^{\mathfrak{r}} a_i\neq 0$. 
\begin{proposition}\label{number of solutions to a quad form}
Let $\ell$ be odd and $Q(x_1,\dots, x_n)\in \F_q[x_1,\dots, x_n]$ be a quadratic form. Setting $\mathfrak{r}:=\op{rank} Q$, we denote by $N(Q=0)$ the number of solutions of $Q=0$. Then, we have that
\[N(Q=0)=\begin{cases}  q^{n-1}+(q-1)q^{n-\mathfrak{r}/2-1}\eta\left((-1)^{\mathfrak{r}/2}\Delta(Q_0)\right),& \text{ if }\mathfrak{r}\text{ is even;}\\
q^{n-1}, & \text{ if }\mathfrak{r}\text{ is odd.}\end{cases}\]
\end{proposition}
\begin{proof}
Let $\mathfrak{r}$ denote the rank of $Q$, and first consider the case when $\mathfrak{r}$ is even. According to Proposition \ref{diagonal form presentation}, $Q$ is equivalent to a diagonal form. Thus, without loss of generality, assume that there are non-zero constants $a_1, \dots, a_{\mathfrak{r}}$ such that
\[Q(x_1, \dots, x_{\mathfrak{r}}, \dots, x_n)=\sum_{i=1}^{\mathfrak{r}} a_i x_i^2.\] Thus, for every solution to 
\[Q_0(x_1, \dots, x_{\mathfrak{r}})=\sum_{i=1}^{\mathfrak{r}} a_i x_i^2=0,\] there are $q^{n-\frak{r}}$ solutions to $Q(x_1, \dots, x_{\mathfrak{r}}, \dots, x_n)=0$. As a consequence, we find that $N(Q=0)=q^{n-\mathfrak{r}} N(Q_0=0)$. Since $Q_0$ is nondegenerate, the result follows from Theorem~\ref{thm:quadratic-Fq}, according to which \[N(Q_0=0)=q^{\mathfrak{r}-1}+(q-1)q^{(\mathfrak{r}-2)/2}\eta\left((-1)^{\mathfrak{r}/2}\Delta(Q_0)\right).\]

On the other hand, if $\mathfrak{r}$ is odd, we have that $N(Q_0=0)=q^{\mathfrak{r}-1}$, according to \cite[Theorem 6.27]{lidl1997finite}.
\end{proof}

\section{Statistics for the Iwasawa invariants of bouquets}\label{S:bouquet}
\par Fix an odd prime number $\ell$. Recall that a \emph{bouquet} is a multigraph with a single vertex. As a multigraph, each undirected edge corresponds to a pair of edges oriented in opposite directions, starting and ending at the unique vertex. Thus in particular, a bouquet is uniquely determined by the number of undirected loops. Let $X_t$ be the bouquet with $t$ undirected loops labelled $e_1,\dots, e_t$. Associated to a vector $\alpha=(\alpha_1,\ldots,\alpha_t)\in\Zl^t$ is a $\Z_\ell$-tower of $X_t$. We say that $\alpha$ is \emph{admissible} if for some coordinate $\alpha_i$, we have that $\ell\nmid \alpha_i$. Note that according to Corollary \ref{admissible corollary for bouquets}, $\alpha$ is admissible if and only if Assumption \ref{main assumption} is satisfied. Let $\tau=\tau_\ell$ be the Haar measure on $\Z_\ell^t$ and note that $\tau\left(\Zl^t\setminus (\ell \Zl)^t\right)=1-\ell^{-t}$.
\par In this section, we study two related questions:
\begin{enumerate}
    \item\label{bouquet question 1} Given a fixed bouquet $X_t$, how do the Iwasawa invariants $\mu$ and $\lambda$ vary, as 
    $\alpha=(\alpha_1,\dots, \alpha_t)$ varies over the set of admissible vectors in $\Z_\ell^t$? It follows from Theorem \ref{lambda is odd} that the $\lambda$-invariant for a bouquet is odd. Thus, in particular, the $\lambda$ invariant is $\geq 1$. We fix constants $\mu_0\in \Z_{\geq 0}$ and $\lambda_0\in \Z_{\geq 1}$, let $\mathcal{S}(\mu_0, \lambda_0,\ell,t)$ be the set of vectors $\alpha\in \Z_\ell^t\backslash (\ell \Z_\ell)^t$ such that $\mu_\ell(X_t,\alpha)=\mu_0$ and $\lambda_\ell(X_t,\alpha)=\lambda_0$.
    This probability is computed as follows
    \[\begin{split}\op{Prob}_{\ell}\left(X_t; \mu=\mu_0, \lambda=\lambda_0\right):= & \frac{\tau\left(\mathcal{S}(\mu_0, \lambda_0,\ell,t)\right)}{\tau\left(\Zl^t\setminus (\ell \Zl)^t\right)}\\
    = & (1-\ell^{-t})^{-1}\tau\left(\mathcal{S}(\mu_0, \lambda_0,\ell,t)\right).\end{split}\] The above formula does make sense since in the cases of interest, we shall implicitly show that $\mathcal{S}(\mu_0, \lambda_0,\ell,t)$ is measurable. Our computations are possible for certain values of $\mu_0$ and $\lambda_0$. The probabilities $\op{Prob}_{\ell}(X_t; \mu>\mu_0)$, $\op{Prob}_{\ell}(X_t; \mu=\mu_0)$, $\op{Prob}_{\ell}(X_t; \mu=\mu_0, \lambda\geq \lambda_0)$ etc. are defined similarly.
    \item\label{bouquet question 2} Let $\mathcal{T}(\mu_0, \lambda_0, \ell)$ be the set of pairs $(t, \alpha)$, where $t\in \Z_{\geq 2}$ and admissible $\alpha=(\alpha_1,\dots, \alpha_t)\in \Z^t$ such that $\mu_\ell(X_t, \alpha)=\mu_0$ and $\lambda_\ell(X_t, \alpha)=\lambda_0$. Fix $\delta\in (0,1)$ and let $\mathcal{T}_{\leq x}(\mu_0, \lambda_0, \ell)$ be the set of $(t,\alpha)\in \mathcal{T}(\mu_0, \lambda_0, \ell)$ such that $t\leq x^\delta$ and $\op{max}\{ |\alpha_1|,\dots, |\alpha_t|\}\leq x$. Let $\mathcal{A}_{\leq x}(\ell)$ be the tuples $(t,\alpha)$ for which $\alpha$ is admissible and such that $t\leq x^\delta$ and $\op{max}\{ |\alpha_1|,\dots, |\alpha_t|\}\leq x$. 
    The density of $\mathcal{T}_{\leq x}(\mu_0, \lambda_0, \ell)$ is given by
    \[\mathfrak{d}_{\ell}(\mu_0, \lambda_0):=\lim_{x\rightarrow \infty} \left(\frac{\# \mathcal{T}_{\leq x}(\mu_0, \lambda_0, \ell)}{\#\mathcal{A}_{\leq x}(\ell)}\right),\] provided the above limit exists. 
\end{enumerate}

\subsection{The vanishing of the $\mu$-invariant}
In this section, we study when the $\mu$-invariant of a $\Zl$-tower of a bouquet $X_t$ vanishes. We are interested in the power series
\begin{equation}\label{def of f(T)}
f(T)=f_{X_t,\alpha}(T)=\sum_{k=1}^tf_{\alpha_k}(T)\in\Zl\lb T\rb,
\end{equation}
where we recall that
\[
f_{a}(T)=2-(1+T)^{a}-(1+T)^{-a}=-\sum_{n\ge2}\left(\binom{a}{n}+\binom{-a}{n}\right)T^n.
\]
Upon expanding the binomial coefficients, we deduce that the coefficient of $T^n$ in $f_a(T)$ is of the form 
\[
\frac{2}{n!}\sum_{i=1}^{\lfloor n/2\rfloor}c_{n,i}a^{2i},
\]
where $c_{n,i}\in\Z$ is independent of $a$. More precisely, $c_{n,i}$ is the coefficient of $x^{2i}$ in $-x(x-1)\dots(x-n+1)$. Furthermore, if $n=2m$ is even, then the highest term in $a$ is $a^{2m}$ with coefficient $c_{2m, m}=- 1$, whereas if $n=2m+1$ is odd, then the highest term is $a^{2m}$ with coefficient $c_{2m+1,m}=m(2m+1)$. We write $f(T)=\sum_{n\geq 2} \beta_{n} T^{n}$, 
in particular, we have
\[
\beta_n=\sum_{k=1}^t\frac{2}{n!}\sum_{i=1}^{\lfloor n/2\rfloor}c_{n,i}\alpha_k^{2i}.
\]

\begin{lemma}\label{mu invariant crit lemma}
Let $\ell$ be an odd prime number and $f(T)$ be the power series associated to the $\mathbb{Z}_{\ell}$-tower of the  bouquet $X_t$ given by the admissible vector $\alpha=(\alpha_1, \dots, \alpha_t)$. Suppose  that the associated $\mu$-invariant $\mu_\ell(X_t, \alpha)$ is positive. Then, for all $m\geq 1$, we have that \[\sum_{k} \alpha_k^{2m}\equiv 0\mod{\ell^{\mu}}.\]
\end{lemma}
\begin{proof}
If $\mu>0$, then, all coefficients $\beta_n$ are divisible by $\ell^{\mu}$. Setting $p_i:=\sum_{k} \alpha_k^{2i}$, we find that 
\[\beta_n=\frac{2}{n!}\sum_{i=1}^{\lfloor n/2\rfloor} c_{n, i} p_i\equiv 0\mod{\ell^{\mu}}.\]Then, we find that $\sum_{i=1}^{\lfloor n/2\rfloor} c_{n, i} p_i\equiv 0\mod{\ell^{\mu}}$. Set $n=2m$, and obtain from the above that 
\[p_m=-c_{2m, m}p_m\equiv \sum_{i=1}^{m-1} c_{2m, i} p_i\mod{\ell^{\mu}}.\]Therefore, by induction on $m$, we find that $p_m\equiv 0\mod{\ell^\mu}$ for all $m$.
\end{proof}

Next, we prove a result on the $\mu$-invariant. Note that since $\alpha$ is assumed to be admissible, we have that $\#\{k:\alpha_k\not\equiv 0\mod \ell\}>0$.

\begin{proposition}\label{ellpowermu}
Assume that $\ell$ is odd. If $\mu>0$, then $\#\{k:\alpha_k\not\equiv 0\mod \ell\}$ is a multiple of $\ell$.
\end{proposition}
\begin{proof}
Let us write $p_i:=\sum_{k=1}^t\alpha_k^{2i}$. By Lemma \ref{mu invariant crit lemma}, we find that $p_i\equiv 0\mod \ell^{\mu}$ for all $i\ge1$.

Let us write $N=\#\{k:\alpha_k\not\equiv 0\mod \ell\}$. For $1\le i\le N$, let $p_i'$ denote the power sum $\displaystyle\sum_{\ell\nmid \alpha_k}\alpha_k^{2i}$ and let $e_i$ denote the symmetric sum $\displaystyle\sum_{k_1<\cdots<k_i} \alpha_{k_1}^2\cdots \alpha_{k_i}^2$ of the same elements. In particular, $p_i'\equiv 0\mod \ell$ for all $i$. Newton's identities (see \cite{mead1992newton}) say that 
\[
ie_i=\sum_{j=1}^i(-1)^{j-1}e_{i-j}p'_j,\quad 1\le i\le N,
\]
where $e_0$ is defined to be $1$.  In particular, $Ne_N\equiv 0\mod \ell$.  By definition, $$ e_N=\prod_{\ell\nmid \alpha_k}\alpha_k^2\not\equiv 0\mod \ell$$ Thus, we deduce that $\ell|N$.
\end{proof}

The result below improves upon \cite[Corollary 5.8]{Vallieres:2021}. 
\begin{corollary}\label{cor:posmu}
Let $\ell$ be an odd prime. If $t<\ell$, then $\mu=0$.
\end{corollary}
\begin{proof}
Note that since there exists some value $\alpha_k$ such that $\ell\nmid \alpha_k$, we find that $\#\{k:\alpha_k\not\equiv 0\mod \ell\}>0$. According to Proposition \ref{ellpowermu}, 
\[t\geq \#\{k:\alpha_k\not\equiv 0\mod \ell\}\geq \ell,\]and the result follows.
\end{proof}

The following result is a refinement of Proposition \ref{ellpowermu}.

\begin{proposition}\label{prop:quadres}
Let $\ell$ be an odd prime, and suppose that $\mu>0$. Let $Q_\ell$ be the set of non-zero quadratic residues modulo $\ell$. For each $x\in Q_\ell$, let  $r_x:=\#\{\alpha_k:\alpha_k^2\equiv x\mod \ell\}$. Then $r_x\equiv0\mod\ell^{\mu}$ for all $x\in Q_\ell$.
\end{proposition}
\begin{proof}

Recall from the proof Proposition~\ref{ellpowermu} that we have
\[
\sum_{k=1}^t \alpha_k^{2i}\equiv0\mod \ell^{\mu},\quad i\ge 1.
\]
This implies that
\[
\sum_{x\in Q_\ell} r_x \cdot x^{i}\equiv 0\mod\ell^{\mu}, \quad i=1,\ldots, (\ell-1)/2.
\]
The determinant of the Vandermonde matrix $(x^i)_{x\in Q_\ell,1\le i\le (\ell-1)/2}$ is nonzero mod $\ell$. Therefore, we deduce that $r_x\equiv 0\mod \ell^{\mu}$ for all $x$.
\end{proof}

\begin{remark}\label{rk:eg-mu}
The converse to the statement in Proposition \ref{prop:quadres} is not true. For instance, let $\ell=t=3$ and $\alpha=(1,1,2)$, we find that
\[\begin{split}f(T)=& 2\left(2-(1+T)-(1+T)^{-1}\right)+\left(2-(1+T)^2-(1+T)^{-2}\right)\\
\equiv & \frac{2T^4}{(T + 1)^2}\mod{3}.\end{split}\]In particular, the $\mu$-invariant is zero. This  also gives a counterexample for the converse of both Lemma \ref{mu invariant crit lemma} and Proposition \ref{ellpowermu}.
\end{remark}
\begin{theorem}\label{thm:mu>0}
Suppose that $\ell$ is odd and let $m=(\ell-1)/2$.
We have
\[
\PP_{\ell}(X_t;\mu>0)\le \sum_{i=1}^{\lfloor t/\ell\rfloor }2^{i\ell}\binom{t}{i\ell}\sum_{a_1+\cdots+a_m=i}\frac{(i\ell)!}{\prod_{i=1}^m(a_j\ell)!(\ell^t-1)},
\]
where the second sum runs over non-negative integers $a_j\ge0$.
\end{theorem}
\begin{proof}
If $\mu>0$, then Proposition~\ref{ellpowermu} tells us that $\#\{k:\alpha_k\not\equiv 0\mod \ell\}$ is a multiple of $\ell$. Since each residue class modulo $\ell$ in $\alpha\in \Zl^t\setminus (\ell\Zl)^t$ has the same measure, it is enough to count the number of elements in $(r_1,\ldots, r_t)\in\Fl^t\setminus\{0\}$ such that $\#\{k:r_k\neq 0\}$ is a multiple of $\ell$ and $\sum_{k=1}^t r_k^{2i}= 0$, for all $i\ge1$. 

For each $1\le i\le \lfloor t/\ell\rfloor$, we may choose $i\ell$ indices $k_1,\dots, k_{i\ell}$ so that $r_j\neq 0$ if and only if $j\in \{k_b\mid b=1,\dots, i\ell\}$. 
By Proposition~\ref{prop:quadres}, amongst these $i\ell$ elements, we have to distribute $r_j^2$ over the $m$ classes in $Q_\ell$ so that the number of times each quadratic residue occurs is a multiple of $\ell$. 

Let $a_1\ell,\ldots,a_{m}\ell$ be the number of times the residue classes in $Q_\ell$ occur. We have
\[
a_1+\cdots+a_{m}=i.
\]
There are $(i\ell)!/\prod_{i=1}^m(a_j\ell)!$ ways to distribute $i\ell$ distinguishable elements in $m$ distinguishable sets so that the number of elements in these sets are given by $a_1\ell,\ldots,a_m\ell$. 
Once we have distributed the $i\ell$ elements amongst the $m$ classes, each element has two possible values modulo $\ell$, giving us $2^{i\ell}$ choices. Summing up all possible choices of $a_j$'s gives the expression on the right-hand side of the inequality.
\end{proof}

\begin{remark}
Our numerical calculations suggest that for a given $\ell$, the right-hand side of the inequality given by Theorem~\ref{thm:mu>0} tends to $\frac{1}{\ell^{(\ell-1)/2}}$ as $t\rightarrow \infty$.

\begin{figure}[H]
\begin{tabular}{ |p{1cm}||p{5cm}|p{5cm}| }
 \hline
 \multicolumn{3}{|c|}{Numerical estimation of the bound in Theorem \ref{thm:mu>0} for $t=1000$} \\
 \hline
 $\ell$&bound& $\ell^{-(\ell-1)/2}$\\
 \hline
 3& 0.333333333333333&$0.\overline{3}$\\
 5&   0.04&0.04\\
 7& 0.00291545189504373&0.00291545189504373\\
 11& $6.20923395131559\times10^{-6}$&$6.20921323059155\times 10^{-6}$\\
 13& $2.07220439470374\times 10^{-7}$&$2.07176211033003\times 10^{-7}$\\
 \hline
\end{tabular}
\end{figure}
\end{remark}

In what follows, we prove a necessary and sufficient condition for $\mu>0$ in the special case where $\alpha_k\in\Z$ for all $k$. Given an integer $a\in\Z$, let us write $d_k(a)$ for the coefficient of $X^k$ in the shifted Chebyshev polynomial $P_a(X)$ {used in \cite[\S2]{mcgownvallieresII}}. More explicitly,
\[
P_a(X)=2-2 T_a(1-X/2)=d_1(a)X+d_2(a)X^2+\cdots+d_a(a)X^a,
\]
where the Chebyshev polynomials $T_a(X)\in\Z[X]$ of the first kind are defined by the relation
\[
T_a(\cos(\theta))=\cos(a\theta)
\]
for all $\theta\in\mathbb{R}$ (see \cite[\S2]{mcgownvallieresII}).
\begin{lemma} \label{unimodular}
Let $m\ge1$ be an integer. Let $b_{1},\ldots,b_{m}$ be distinct positive integers, then 
$$(d_{b_{i}}(b_{j})) \in {\rm GL}(m,\mathbb{Z}). $$
\end{lemma}
\begin{proof}
Indeed, if we list the $b_{i}$ in increasing order, then the matrix $(d_{b_{i}}(b_{j}))$ is upper-triangular with $\pm 1$ on the diagonal, since $d_{n}(n) = (-1)^{n-1}$ for all $n \in \mathbb{N}$.
\end{proof}

\begin{proposition}\label{prop:pre-image}
Let $\alpha:S \to \mathbb{Z}$ be a function for which there exists at least one $s \in S$ such that $(\alpha(s),\ell)=1$.  Let $b_{1},\ldots,b_{m}$ be the nonzero distinct elements of the set $\{|\alpha(s)|:s \in S\}$.  For $j=1,\ldots,m$, let $t_{j} = \#\{s \in S : |\alpha(s)|=b_{j} \}$.  Then, the polynomial $\sum_{s\in S}P_{\alpha(s)}(X)$ belongs to $\ell\Z[X]$ if and only if $\ell \, | \, t_{j}$ for all $j=1,\ldots,m$. 
\end{proposition}
\begin{proof}
We have $P_{a}(X) = P_{-a}(X)$ for all $a \in \mathbb{Z}$.  Thus, 
\begin{equation} \label{key_eq}
\sum_{s\in S}P_{\alpha(s)}(X) = \sum_{j=1}^{m}t_{j}P_{b_{j}}(X).
\end{equation}
It is clear that $\ell \, | \, t_{j}$ for all $j=1,\ldots,m$ implies that the polynomial above belongs to  $ \ell \mathbb{Z}[X]$.

Conversely, if  $\sum_{s\in S}P_{\alpha(s)}(X) \in \ell \mathbb{Z}[X]$, then \eqref{key_eq} implies that
\begin{equation} \label{syst}
\sum_{j=1}^{m}d_{k}(b_{j})t_{j} \equiv 0 \pmod{\ell}, 
\end{equation}
for all $k \ge 1$.  If we let $b$ be the maximal value of the $b_{j}$, then the matrix $(d_{k}(b_{j}) + \ell\mathbb{Z}) \in M_{b\times m}(\mathbb{F}_{\ell})$ has rank $m$ over $\mathbb{F}_{\ell}$ by Lemma \ref{unimodular}.  Therefore, the system of linear equations \eqref{syst} has a unique solution modulo $\ell$.  This tells us that $\ell \, | \, t_{j}$ for all $j=1,\ldots,m$ and concludes the proof.
\end{proof}

\begin{corollary}
Let $b_1,\ldots,b_m$ be the nonzero distinct elements of the set $\{|\alpha(s)|:s\in S\}$ and write $t_j=\#\{s\in S:|\alpha(s)|=b_j\}$ for $j=1,\ldots, m$. Then $\mu_\ell(X_t,\alpha)>0$ if and only if $\ell|t_j$ for all $j=1,\ldots,m$.
\end{corollary}
\begin{proof}
By \cite[Theorem~5.6]{Vallieres:2021}, $\mu_\ell(X_t,\alpha)>0$ if and only if  $\sum_{s\in S}P_{\alpha(s)}(X)\in\ell\Z[X]$. Therefore, the result follows from Proposition~\ref{prop:pre-image}.
\end{proof}

\subsection{The distribution of Iwasawa invariants of a fixed $X_t$ varying $\alpha$}
\par In this subsection, we fix an odd prime number $\ell$ and fix an integer $t\geq 2$. We let $X_t$ be the bouquet with $t$ undirected loops and we  study question \eqref{bouquet question 1} above. Note that Assumption \ref{euler char assumption} requires that $\chi(X_{t})\neq 0$, and hence we avoid the case when $t=1$; see \S\ref{S:smallt} below for a separate discussion on this case.

We begin by considering $\mathcal{S}(\mu_0, \lambda_0, \ell, t)$. Throughout, it shall be assumed that $\alpha:S\rightarrow \Z_\ell$ is an admissible function. We compute the probability $\op{Prob}(X_t;\mu=\mu_0,\lambda= \lambda_0)$ for various values of $(\mu_0, \lambda_0)\in \Z_{\geq 0}\times \Z_{\geq 0}$. Note that according to Theorem \ref{lambda is odd}, $\lambda$ is odd and thus, the minimal values of $\mu$ and $\lambda$ are $0$ and $1$, respectively.

\begin{lemma}\label{mu=0 lambda=n}
Suppose that $\ell$ is a prime number and $f(T)$ is the power series defined above, see \eqref{def of f(T)}. Then, $\lambda=\lambda_{\ell}(X_t, \alpha)\geq 1$ and the following conditions are equivalent
\begin{enumerate}
    \item $\mu=0$ and $\lambda=n$; 
    \item $\ell\mid \beta_i$ for all $i \leq n$ and $\beta_{n+1}\in \Z_\ell^\times$.
\end{enumerate}
\end{lemma}

\begin{proof}
\par Theorem \ref{lambda is odd} asserts that $\lambda$ is odd, thus $\lambda\geq 1$. Alternatively, this may be directly shown as follows. Express $f(T)$ as $\ell^{\mu} \mathcal{P}(T) u(T)$, where $\mathcal{P}(T)$ is a distinguished polynomial and $u(T)$ is a unit in $\Lambda$. Here $\mu=\mu_{\ell}(X_t, a)$ is the $\mu$-invariant and $\lambda=\lambda_{\ell}(X_t, a)=\op{deg} \mathcal{P}-1$. Note that since $\beta_i=0$ for $i\leq 2$, it follows that $T^2$ divides $\mathcal{P}(T)$, and thus $\op{deg} \mathcal{P}(T)\geq 2$ and $\lambda\geq 1$. 
\par Since $T^2$ divides $\mathcal{P}(T)$, we have that $\mu=0$ and $\lambda=1$, if and only if $f(T)=T^2u(T)$, where $u(T)$ is a unit in $\Lambda$. This is the case precisely when $\beta_2$ is a unit in $\Z_{\ell}$. More generally, if $\mu=0$ and $\lambda=n$, then, $f(T)=\mathcal{P}(T) u(T)$. Reducing modulo $\ell$, we find that $\mathcal{P}(T)\equiv T^{n+1}\mod \ell$, whereby, $\ell\mid \beta_i$ for $i\leq n$ and $\beta_{n+1}\in \Z_\ell^{\times}$. The converse follows from the same argument.
\end{proof}

\begin{lemma}\label{big lemma}
Let $\ell$ be an odd prime number, set $\mu=\mu_\ell(X_t, \alpha)$ and $\lambda=\lambda_\ell(X_t, 
\alpha)$. Then, the following assertions hold.
\begin{enumerate}
    \item\label{lemma 4.3 ass 1} Suppose that $k$ is an integer such that $0<2k<\ell$. Then, the following assertions are equivalent:
    \begin{enumerate}
        \item\label{Lemma c1} $\mu=0$ and $\lambda=2k-1$,
        \item\label{Lemma c2} for all $i<k$, 
        \[\sum_{j=1}^t \alpha_j^{2i}\equiv 0\mod{\ell}\]
        and 
         \[\sum_{j=1}^t \alpha_j^{2k}\not\equiv 0\mod{\ell}.\]
    \end{enumerate}
    \item\label{lemma 4.3 ass 2} If $2k\geq \ell$,  $\mu=0$ and $\lambda=2k-1$, then for all $i<k$, 
        \[\sum_{j=1}^t \alpha_j^{2i}\equiv 0\mod{\ell}.\]
\end{enumerate}
\end{lemma}

\begin{proof}
First, assume that $\mu=0$ and $\lambda=2k-1$. Setting $p_i:=\sum_{j} \alpha_j^{2i}$, we find that 
\[\beta_n=\frac{2}{n!}\sum_{i=1}^{\lfloor n/2\rfloor} c_{n, i} p_i.\] According to Lemma \ref{mu=0 lambda=n}, $\beta_n$ is divisible by $\ell$ for all values of $n<2k$, and $\beta_{2k}$ is not divisible by $\ell$. The assertion of part \eqref{lemma 4.3 ass 2} follows immediately from this. 

We now prove \eqref{lemma 4.3 ass 1}. In particular, we assume that $2k<\ell$. Then, $\frac{2}{n!}$ is an $\ell$-adic unit for all integers $n$ such that $0\le n\leq 2k$. 

Suppose that (b) holds. Then,  $p_i$ is divisible by $\ell$ for all $i<k$, whereas $p_k$ is not divisible by $\ell$. It follows that $\ell$ divides $\beta_n$ for all  $n<2k$, but $\beta_{2k}$ is not divisible by $\ell$. Thus, (a) holds by Lemma \ref{mu=0 lambda=n}.

Conversely, suppose that (a) holds. Then 
\[    \ell\ | \sum_{i=1}^{\lfloor n/2\rfloor}c_{n,i}p_i,\ n<2k,\quad    \ell \nmid \sum_{i=1}^{k}c_{2k,i}p_i.\]
Recall that $c_{2m,m}=-1$ for all $m\ge1$. An inductive argument shows that $\ell|p_i$ for $i<k$, and $\ell\nmid p_k$. In particular, (b) holds.
\end{proof}

\begin{theorem}\label{thm:small-t}
We let $\ell$ be an odd prime and $\mu$ and $\lambda$ denote the Iwasawa invariants $\mu_{\ell}(X_t, \alpha)$ and $\lambda_{\ell}(X_t, \alpha)$, respectively. If $\ell>t$, then $\lambda< 2t$.
\end{theorem}
\begin{proof}By Corollary~\ref{cor:posmu}, $\ell>t$ implies that the $\mu$-invariant of $f(T)$ is zero. Assume by way of contradiction that $\lambda\geq 2t$. Then, $\beta_2,\ldots, \beta_{2t}\equiv 0\mod \ell$. Under the same notation as in the proof of Proposition~\ref{ellpowermu}, we have
\[
p_i\equiv 0\mod \ell, \quad i=1,2,\ldots, t
\]
forcing $Ne_N\equiv 0\mod\ell$ once again by Newton's identities  (see \cite{mead1992newton}).
But this is impossible if $\ell>t$, proving the  affirmation of the theorem.
\end{proof}
\begin{remark}\label{rk:small-egs}
Without the condition $\ell>t$, it is possible that $\mu=0$ and $\lambda> 2t$.  In fact, applying Lemma \ref{arb_large} with $n_{1}=0$ shows that it is possible to find examples with $\mu=0$ and $\lambda$ arbitrarily large when $t = \ell+1$.  When $t=\ell=2$, we have found examples with $\mu=0$ and $$\lambda = 5,9,17,33,65,129,257.$$ 
Whereas when $t=3$, for $\ell=2$ we have found examples with $\mu=0$ and 
$$\lambda=3,5,7,9,15,17,31,33,63,65,127,129,255,257,$$
and for $\ell = 3$ we have found examples with $\mu=0$ and 
$$\lambda = 3,5,9,17,27,53,81,161,243.$$
For instance, if we take $t=3$, $\ell = 3$, with the voltage assignment on the three edges given by 1, 8 and 10, respectively, then we get:
\begin{equation*}
\begin{tikzpicture}[baseline={([yshift=-1.7ex] current bounding box.center)}]
\node[draw=none,minimum size=3cm,regular polygon,regular polygon sides=1] (a) {};
\foreach \x in {1}
  \fill (a.corner \x) circle[radius=0.7pt];
\draw (a.corner 1) to [in=50,out=130,loop] (a.corner 1);
\draw (a.corner 1) to [in=50,out=130,distance = 0.8cm,loop] (a.corner 1);
\draw (a.corner 1) to [in=50,out=130,distance = 0.5cm,loop] (a.corner 1);
\end{tikzpicture}
\longleftarrow \, \, \,
\begin{tikzpicture}[baseline={([yshift=-0.6ex] current bounding box.center)}]
\node[draw=none,minimum size=2cm,regular polygon,regular polygon sides=3] (a) {};

\foreach \x in {1,2,3}
  \fill (a.corner \x) circle[radius=0.7pt];

\path (a.corner 1) edge [bend left=20] (a.corner 2);
\path (a.corner 1) edge [bend right=20] (a.corner 2);
\path (a.corner 2) edge [bend left=20] (a.corner 3);
\path (a.corner 2) edge [bend right=20] (a.corner 3);
\path (a.corner 3) edge [bend left=20] (a.corner 1);
\path (a.corner 3) edge [bend right=20] (a.corner 1);

\path (a.corner 1) edge  (a.corner 2);
\path (a.corner 2) edge  (a.corner 3);
\path (a.corner 3) edge  (a.corner 1);

\end{tikzpicture}
\longleftarrow \, \,
\begin{tikzpicture}[baseline={([yshift=-0.6ex] current bounding box.center)}]
\node[draw=none,minimum size=2cm,regular polygon,regular polygon sides=9] (a) {};

\foreach \x in {1,2,...,9}
  \fill (a.corner \x) circle[radius=0.7pt];
  
\path (a.corner 1) edge [bend left=20] (a.corner 2);
\path (a.corner 1) edge [bend right=20] (a.corner 2);
\path (a.corner 2) edge [bend left=20] (a.corner 3);
\path (a.corner 2) edge [bend right=20] (a.corner 3);
\path (a.corner 3) edge [bend left=20] (a.corner 4);
\path (a.corner 3) edge [bend right=20] (a.corner 4);
\path (a.corner 4) edge [bend left=20] (a.corner 5);
\path (a.corner 4) edge [bend right=20] (a.corner 5);
\path (a.corner 5) edge [bend left=20] (a.corner 6);
\path (a.corner 5) edge [bend right=20] (a.corner 6);
\path (a.corner 6) edge [bend left=20] (a.corner 7);
\path (a.corner 6) edge [bend right=20] (a.corner 7);
\path (a.corner 7) edge [bend left=20] (a.corner 8);
\path (a.corner 7) edge [bend right=20] (a.corner 8);
\path (a.corner 8) edge [bend left=20] (a.corner 9);
\path (a.corner 8) edge [bend right=20] (a.corner 9);
\path (a.corner 9) edge [bend left=20] (a.corner 1);
\path (a.corner 9) edge [bend right=20] (a.corner 1);

\path (a.corner 1) edge  (a.corner 2);
\path (a.corner 2) edge  (a.corner 3);
\path (a.corner 3) edge  (a.corner 4);
\path (a.corner 4) edge  (a.corner 5);
\path (a.corner 5) edge  (a.corner 6);
\path (a.corner 6) edge  (a.corner 7);
\path (a.corner 7) edge  (a.corner 8);
\path (a.corner 8) edge  (a.corner 9);
\path (a.corner 9) edge  (a.corner 1);

\end{tikzpicture}
\longleftarrow
\begin{tikzpicture}[baseline={([yshift=-0.6ex] current bounding box.center)}]
\node[draw=none,minimum size=2cm,regular polygon,regular polygon sides=27] (a) {};

\foreach \x in {1,2,...,27}
  \fill (a.corner \x) circle[radius=0.7pt];
  
\foreach \y\z in {1/2,2/3,3/4,4/5,5/6,6/7,7/8,8/9,9/10,10/11,11/12,12/13,13/14,14/15,15/16,16/17,17/18,18/19,19/20,20/21,21/22,22/23,23/24,24/25,25/26,26/27,27/1}
  \path (a.corner \y) edge (a.corner \z);
  
\foreach \y\z in {1/9,2/10,3/11,4/12,5/13,6/14,7/15,8/16,9/17,10/18,11/19,12/20,13/21,14/22,15/23,16/24,17/25,18/26,19/27,20/1,21/2,22/3,23/4,24/5,25/6,26/7,27/8}
  \path (a.corner \y) edge (a.corner \z); 
  
\foreach \y\z in {1/11,2/12,3/13,4/14,5/15,6/16,7/17,8/18,9/19,10/20,11/21,12/22,13/23,14/24,15/25,16/26,17/27,18/1,19/2,20/3,21/4,22/5,23/6,24/7,25/8,26/9,27/10}
  \path (a.corner \y) edge (a.corner \z);

\end{tikzpicture}
\, \, \longleftarrow \ldots
\end{equation*}
The power series $f(T)$ starts as follows 
$$f(T) = -165T^2 +165T^{3}-1326T^{4} + \ldots +3470655T^{17} - 5167526T^{18}+ \ldots,$$
where $3$ divides all the coefficients up to and including $3470655$, but it does not divide $5167526$. Thus, $\mu = 0$ and $\lambda = 17$.  Furthermore, since
$$\log_{3}\left(\frac{3}{2} \cdot 18 \right) = 3,$$
\cite[Theorem 4.1]{mcgownvallieresII} implies that if $n \ge 3$, we have
$${\rm ord}_{3}(\kappa_{n}) = 17n + \nu $$
for some $\nu \in \mathbb{Z}$.  Using SageMath \cite{SAGE}, we calculate
$$\kappa_{0} = 1, \kappa_{1} = 3^{3} , \kappa_{2} =  3^{10}, \kappa_{3} = 2^{18} \cdot 3^{27}, $$
and
$$\kappa_{4} = 2^{18} \cdot 3^{44} \cdot 163^{2} \cdot 487^{2} \cdot 37907^{2} \cdot 799471^{2}. $$
Thus, we have
$${\rm ord}_{3}(\kappa_{n}) =  17n - 24,$$
for all $n \ge 2$.
\end{remark}
\begin{remark}
It directly follows from Theorem~\ref{thm:small-t} that if $\ell>t$, then $\op{Prob}_{\ell}(X_t; \mu=\mu_0, \lambda=\lambda_0)=0$ unless $\mu_0=0$ and $\lambda_0< 2t$.
\end{remark}

We wish to estimate the probability $\op{Prob}_\ell(X_t; \mu=0, \lambda=\lambda_0)$. Note that by Theorem~\ref{lambda is odd}, $\op{Prob}_\ell(X_t; \lambda=\lambda_0)=0$ whenever $\lambda_0$ is even. Therefore, we assume that $\lambda_0=2k-1$ is odd.
\begin{theorem}\label{thm:mu0lambda small}
Let $X_t$ be a bouquet and $k>1$. Assume that $\ell$ is odd, $2k-1<\ell$ and $k(k-1)<t$. Then, we have that 
\[\op{Prob}_\ell(X_t; \mu=0, \lambda< 2k-1)\leq (1-\ell^{-t})^{-1}(1-\ell^{-k(k-1)}).\]

\end{theorem} 
\begin{proof}
Note that since $\ell$ is assumed to be odd and $2k-1<\ell$, it is automatically true that $2k-1<\ell-1$. Let $\alpha:S\rightarrow \Z_\ell$ be an admissible function for which $\mu=0$ and $\lambda<2k-1$. We set $\lambda_0:=2k-1$, note that by assumption, $\lambda_0<\ell-1$. Let $\bar{\alpha}=(\bar{\alpha}_1, \dots, \bar{\alpha}_t)$ be the mod-$\ell$ reduction of $\alpha$. According to Lemma \ref{big lemma}, $\bar{\alpha}$ is \emph{not} a solution to all of the equations
\[Q_i(x_1,\dots, x_t)=\sum_{j=1}^{t} x_j^{2i}= 0\]for $i<k$. Then, since $\lambda_0<\ell-1$, it follows that $\beta_n$ is divisible by $\ell$ for $n\leq \lambda_0$. Let $N$ be the number of non-zero solutions $\bar{\alpha}$ such that $Q_i(\bar{\alpha})=0$ for all $i<k$. We find that 
\[\op{Prob}_\ell(X_t; \mu=0, \lambda< 2k-1)\leq 1-\frac{N}{\ell^{t}-1}.\]
The sum of degrees $\sum_{i=1}^{k-1} \op{deg} (Q_i)=\sum_{i=1}^{k-1} 2i=k(k-1)$. Since $k(k-1)<t$ by assumption, it follows from Theorem \ref{lower bound on number of solns} that $N\geq \ell^{t-k(k-1)}-1$, and therefore,
\[\op{Prob}_\ell(X_t; \mu=0, \lambda< 2k-1)\leq 1-\frac{\ell^{t-k(k-1)}-1}{\ell^{t}-1}=(1-\ell^{-t})^{-1}(1-\ell^{-k(k-1)}).\]

\end{proof}

\begin{example}Let us now illustrate the above result through an example. Let $t=10$ and  $k=3$. We study the proportion of admissible vectors for which $\mu=0$ and $\lambda<5$. Note that $k(k-1)=6<t=10$. Thus, for all primes $\ell>2k-1=7$, we have that \[\op{Prob}_\ell(X_t; \mu=0, \lambda< 5)\leq (1-\ell^{-10})^{-1}(1-\ell^{-6}).\]
\end{example}
\par Next, we study the proportion of admissible vectors 
$\alpha=(\alpha_1,\dots, \alpha_t)$ such that $\mu=0$ and $\lambda=1$. This case is of particular interest, since it is indeed the case when the Iwasawa invariants take on their minimal values.  

\begin{theorem}\label{thm:mu0lambda1}
Let $\ell$ be an odd prime number and $t \in \mathbb{Z}_{\ge 2}$. Then, we have that 
\[\op{Prob}_{\ell}(X_t;\mu=0, \lambda=1)=\begin{cases} 1-\frac{\ell^{t-1}+(-1)^{\frac{t(\ell-1)}{4}}(\ell-1)\ell^{\frac{t}{2}-1}-1}{\ell^{t}-1},& \text{ if }t \text{ is even};
\\ 
1-\frac{\ell^{t-1}-1}{\ell^t-1} & \text{ if }t\text{ is odd.}\end{cases}\]
\end{theorem}
\begin{proof}
Consider the quadratic form $Q(x_1,\dots, x_t)=\sum_{i=1}^t x_i^2$ defined over $\F_\ell$. Then $\mu=0$ and $\lambda=1$ if and only if $Q(\bar a_1,\ldots,\bar a_t)\ne0$. Therefore,
\[
\op{Prob}_{\ell}(X_t;\mu=0, \lambda=1)=1-\frac{N(Q=0)-1}{\ell^t-1}.
\]
Therefore, the theorem is now a direct consequence of Proposition \ref{number of solutions to a quad form}.
\end{proof}

\subsection{Results for small values of $t$}\label{S:smallt}
\par In this subsection, we study the distribution of the $\lambda$-invariant in greater detail for various values of $t$, starting with $t=1$. 
Throughout, it is assumed that $\alpha=(\alpha_1, \dots, \alpha_t)$ is admissible.
\subsubsection*{$\mathbf{t=1}$} 
The bouquet $X_1$ consists of a single loop. Let $\ell$ be any prime number ($\ell=2$ is allowed). Even though  Assumption \ref{euler char assumption} does not hold,  we can compute  \[f_a(T)=\beta_2 T^2+\beta_3 T^3+\beta_4 T^4+\dots,\] to find that $\beta_2=\alpha^2$. Since $\ell\nmid \alpha$, we find that $\beta_2\in \Z_\ell^\times$, and thus, $\mu=0$ and $\lambda=1$ for all choices of $a\in \Z_\ell^{\times}$. In particular, we find for \emph{all} primes $\ell$, 
\[\op{Prob}_{\ell}(X_1;\mu, \lambda)=\begin{cases}
1&\text{ if }\mu=0\text{ and }\lambda=1,\\
0&\text{ otherwise.}
\end{cases}\] 
The Iwasawa invariants are seen to be $\mu=0$ and $\lambda=1$ for all choices of $\alpha$.
\subsubsection*{$\mathbf{t=2}$}
  Let $\alpha=(\alpha_1,\alpha_2)\in \Z_\ell^2$ be an admissible vector. Suppose that $\ell\geq 3$. It follows from Theorems~\ref{lambda is odd} and \ref{thm:small-t} that $\mu=\mu_{\ell}(X_2, \alpha)=0$ and $\lambda=\lambda_{\ell}(X_2, \alpha)\in\{1,3\}$. Combined with Theorem~\ref{thm:mu0lambda1}, we deduce the following.

\begin{proposition}
Suppose that $\ell\ge3$ and $t=2$. Then, we have that $\mu=\mu_{\ell}(X_2, \alpha)=0$ for all choices of $\alpha\in \Z_\ell^2\backslash (\ell\Z_\ell)^2$. If $\ell \equiv 3\mod{4}$, then we always have $\lambda=1$. If $\ell \equiv 1\mod{4}$, we have \[\op{Prob}_{\ell}(X_2; \mu=0,\lambda=\lambda_0)=\begin{cases}
    \frac{\ell-1}{\ell+1}&\text{ if }\lambda_0=1,\\
    \frac{2}{\ell+1}&\text{ if }\lambda_0=3,\\
    0&\text{ otherwise.}\\
    \end{cases}\]
\end{proposition}

\subsubsection*{$\mathbf{t=3}$} We now suppose that $\ell\ge5$ (the case $\ell=t=3$ seems a lot more subtle; see Remark~\ref{rk:eg-mu} and \ref{rk:small-egs}). As before, Theorem~\ref{thm:small-t} tells us that the $\mu$-invariant is always zero. We study the distribution of $\lambda$-invariants in the following proposition.

\begin{proposition}
Let $\ell \ge 5$ be a prime number and let $t=3$. Then, we have that 
\[\op{Prob}_{\ell}(X_3; \mu=0,\lambda=1)=\frac{\ell^2}{\ell^2+\ell+1}.\]
Suppose that $\ell\equiv 2\mod{3}$, then we have that 
\[\op{Prob}_{\ell}(X_3; \mu=\mu_0, \lambda=\lambda_0)=\begin{cases}
\frac{\ell^2}{\ell^2+\ell+1}&\text{ if }\mu_0=0\text{ and }\lambda_0=1,\\
\frac{\ell+1}{\ell^2+\ell+1}&\text{ if }\mu_0=0\text{ and }\lambda_0=3,\\
0&\text{ otherwise.}
\end{cases}\]
Suppose that $\ell \equiv 1\mod{3}$, then we have that
\[\op{Prob}_{\ell}(X_3; \mu=\mu_0, \lambda=\lambda_0)=\begin{cases}
\frac{\ell^2}{\ell^2+\ell+1}&\text{ if }\mu_0=0\text{ and }\lambda_0=1,\\
\frac{\ell-7}{\ell^2+\ell+1}&\text{ if }\mu_0=0\text{ and }\lambda_0=3,\\
\frac{8}{\ell^2+\ell+1}&\text{ if }\mu_0=0\text{ and }\lambda_0=5,\\
0&\text{ otherwise.}
\end{cases}\]
\end{proposition}
\begin{proof}
Recall from Theorem~\ref{thm:mu0lambda1} that the probability of $\mu=0$ and $\lambda=1$ is given by
\[
\frac{\ell^2}{\ell^2+\ell+1}.
\]
Furthermore, Theorems~\ref{lambda is odd} and \ref{thm:small-t} imply that $\lambda\in\{1,3,5\}$. More specifically, we have
\begin{align*}
\lambda=1&\Leftrightarrow  \alpha_1^2+ \alpha_2^2+\alpha_3^2\not\equiv 0\mod \ell;\\    
\lambda=3&\Leftrightarrow \alpha_1^2+\alpha_2^2+ \alpha_3^2\equiv 0\mod \ell,\quad \alpha_1^4+\alpha_2^4+\alpha_3^4\not\equiv 0\mod \ell ;\\
\lambda=5&\Leftrightarrow \alpha_1^2+\alpha_2^2+\alpha_3^2\equiv 0\mod \ell,\quad \alpha_1^4+ \alpha_2^4+ \alpha_3^4\equiv 0\mod \ell .
\end{align*}

Suppose that $\alpha_1^2+\alpha_2^2+\alpha_3^2\equiv 0\mod \ell$ and that $\alpha_1^4+\alpha_2^4+\alpha_3^4\equiv 0\mod \ell$. By Newton's identities  (see \cite{mead1992newton}), we can replace the second equation by
\[
\alpha_1^2\alpha_2^2+\alpha_2^2\alpha_3^2+\alpha_3^2\alpha_1^2\equiv0\mod\ell.
\]
This gives
\[
\alpha_1^2(\alpha_2^2+\alpha_3^2)+\alpha_2^2\alpha_3^2\equiv \alpha_1^2(-\alpha_1^2)+\alpha_2^2\alpha_3^2\equiv 0\mod \ell.
\]
In particular, we have $\alpha_1^2\equiv \pm \alpha_2\alpha_3$ and so we are reduced to solving 
\begin{equation}
    \alpha_2^2\pm \alpha_2\alpha_3 +\alpha_3^2\equiv 0\mod \ell.\label{eq:a2a3}
\end{equation}

Note that if $\alpha_3\equiv0$, then both $\alpha_1$ and $\alpha_2$ would be divisible by $\ell$, which is not allowed. So, we can divide out by $\alpha_3^{2}$ and we are now reduced to solving
\begin{align}
    X^2\pm X+1&\equiv 0\mod\ell\notag\\
    (2X\pm1)^2&\equiv -3\mod \ell,\label{eq:-3square}
\end{align}
where $X=\alpha_2/\alpha_3$.
In particular, if $\left(\frac{-3}{\ell}\right)=-1$ ($\Leftrightarrow \ell\equiv 2\mod 3$), \eqref{eq:-3square} has no solution.  Thus, it is impossible for $\lambda$ to be $5$. The probability of $\mu=0$ and $\lambda=3$ is therefore
\[
1-\frac{\ell^2}{\ell^2+\ell+1}=\frac{\ell+1}{\ell^2+\ell+1}.
\]

If $\left(\frac{-3}{\ell}\right)=+1$ ($\Leftrightarrow \ell\equiv 1\mod 3$), then \eqref{eq:-3square}  admits two solutions. Each choice of $X$ gives rise to $\ell-1$ pairs of $(\bar \alpha_2,\bar \alpha_3)\in\Fl^2$. Taking into account the two choices of signs, we have in total $4(\ell-1)$ choices of $(\bar \alpha_2,\bar \alpha_3)$ satisfying \eqref{eq:a2a3}. It remains to choose $\alpha_1$ such that
\[
\alpha_1^2\equiv\pm \alpha_2\alpha_3\equiv\pm X\alpha_3^2.
\]
Note that $\pm X$ is a third root of unity, which is always a square in $\Fl$ (since $\ell\equiv 1\mod 6$, implying that $\Fl$ admits a sixth root of unity). Thus, we have two choices of $\bar \alpha_1$ for each pair of $(\bar \alpha_2,\bar \alpha_3)$.
In total, we have $8(\ell-1)$ triples of $(\bar \alpha_1,\bar \alpha_2,\bar \alpha_3)$ satisfying $\bar \alpha_1^2+\bar \alpha_2^2+\bar \alpha_3^2=\bar \alpha_1^4+\bar \alpha_2^4+\bar \alpha_3^4=0$.
Therefore, the probability $\lambda=5$ is
\[
\frac{8(\ell-1)}{\ell^3-1}=\frac{8}{\ell^2+\ell+1},
\]
whereas that of $\mu=0$ and $\lambda=3$ is 
\[
1-\frac{\ell^2}{\ell^2+\ell+1}-\frac{8}{\ell^2+\ell+1}=\frac{\ell-7}{\ell^2+\ell+1}.
\]
\end{proof}

\subsection{Vary $t$ and $\alpha$}

\par In this subsection, we study question \eqref{bouquet question 2} mentioned at the start of the section. In other words, we prove results about the variation of $\mu_{t,\alpha}=\mu_\ell(X_t,\alpha)$ and $\lambda_{t,\alpha}=\lambda_\ell(X_t,\alpha)$ when both $t$ and $\alpha$ vary. Recall that the prime $\ell$ is fixed. We shall assume that $\ell$ is odd. Given $t\in \Z_{\geq 2}$, define $\mathcal{T}(\mu,\lambda,\ell,t)$ to be the set of integral vectors $\alpha = (\alpha_1,\ldots,\alpha_t)\in \Z^t$ such that
\begin{itemize}
        \item for some coordinate $\alpha_i$, we have $\ell\nmid \alpha_i$, 
        \item $\mu_{t,\alpha}=\mu$ and $\lambda_{t,\alpha}=\lambda$.
\end{itemize}
Set $\mathcal{T}_{\leq x}(\mu,\lambda,\ell,t)$ to be the set of $\alpha \in \mathcal{T}(\mu,\lambda,\ell,t)$ such that $\max\{|\alpha_1|,\cdots, |\alpha_t|\} \leq x$. Let $\mathcal{A}_{\leq x}(\ell, t)$ be the set of all tuples $\alpha$ such that
\begin{itemize}
    \item $\max\{|\alpha_1|,\cdots, |\alpha_t|\}\leq x$,
    \item $\ell\nmid \alpha_i$ for some $1\leq i\leq t$. 
\end{itemize} 

On the other hand, fix $\delta\in (0,1)$ and set 
\[\mathcal{T}_{\leq x}(\mu, \lambda, \ell):=\bigcup_{t\leq x^\delta} \mathcal{T}_{\leq x}(\mu, \lambda, \ell,t)\text{ and } \mathcal{A}_{\leq x}(\ell):=\bigcup_{t\leq x^{\delta}} \mathcal{A}_{\leq x}(\ell,t).\] We now state the main result of this section.
\begin{theorem}\label{section 3 main thm}
With respect to notation above, 
\[\lim_{x\rightarrow \infty}\frac{\# \mathcal{T}_{\leq x}(0, 1, \ell)}{\# \mathcal{A}_{\leq x}(\ell)}=1-\ell^{-1}. \]
\end{theorem}

\begin{remark}
We note here that the limit above is equal to the limit of the densities in Theorem \ref{thm:mu0lambda1} as $t\rightarrow \infty$. This does not follow from formal arguments and further work is required to obtain the result. In fact, it is crucial that the quantity $\delta$ defined above lies in the open interval $(0,1)$. The method used in proving Theorem \ref{section 3 main thm} does not extend to $\delta\geq 1$. 
\end{remark}

\begin{lemma}\label{lem:bound_A}
With respect to notation above,
\[\left|\#\mathcal{A}_{\leq x}(\ell, t)-2^t\left(1-\ell^{-t}\right)x^t\right|< ctx^{t-1}\] for some absolute constant $c>0$ which is independent of $t$ and $x$.
\end{lemma}
\begin{proof}
The number of $\alpha$ such that $\op{max}\left\{|\alpha_1|, |\alpha_2|, \dots, |\alpha_t|\right\}\leq x$ is $\left(2\lfloor x\rfloor+1\right)^t$. The number of tuples such that such that $\op{max}\left\{|\alpha_1|, |\alpha_2|, \dots, |\alpha_t|\right\}\leq x$ and $\ell$ divides all $\alpha_i$ is $\left(2\lfloor x/\ell\rfloor+1\right)^t$. Therefore, we find that \[\#\mathcal{A}_{\leq x}(\ell, t)=\left(2\lfloor x\rfloor+1\right)^t-\left(2\lfloor x/\ell\rfloor+1\right)^t.\]

We have that 
\[\begin{split}
    & \left|\#\mathcal{A}_{\leq x}(\ell, t)-2^t\left(1-\ell^{-t}\right)x^t\right|\\
    = & \left| \left(2\lfloor x\rfloor+1\right)^t-\left(2\lfloor x/\ell\rfloor+1\right)^t-2^t\left(1-\ell^{-t}\right)x^t\right|\\
    \leq & \left|(2\lfloor x\rfloor+1)^t-(2x)^t\right|+\left|(2\lfloor x/\ell\rfloor+1)^t-(2x/\ell)^t\right|\\
    \leq & t\left|(2\lfloor x\rfloor+1)\right|^{t-1}+t\left|(2\lfloor x/\ell\rfloor+1)\right|^{t-1}<ct x^{t-1}\\
\end{split}\]
for an absolute constant $c>0$.
In the above inequalities, we have used the mean value theorem to obtain the bound 
\[\left|(2\lfloor x\rfloor+1)^t-(2x)^t\right|\leq t\left|(2\lfloor x\rfloor+1)\right|^{t-1}\]
and similarly with $x$ replaced by $x/\ell$.
\end{proof}

\begin{proposition}\label{prop 4.20}
With respect to the notation above, there is an absolute constant $c>1$ such that
\[\left| \# \mathcal{A}_{\leq x}(\ell) - (1-\ell^{-\lfloor x^\delta\rfloor})(2x)^{\lfloor x^\delta\rfloor}\right|< c (2x)^{\lfloor x^\delta\rfloor -1}.\]
\end{proposition}
\begin{proof}
Note that 
\[\# \mathcal{A}_{\leq x}(\ell)=\sum_{t=1}^{\lfloor x^\delta\rfloor} \# \mathcal{A}_{\leq x}(\ell,t).\]
Therefore, we deduce from Lemma~\ref{lem:bound_A} that 
\[\begin{split}
 &  \left| \# \mathcal{A}_{\leq x}(\ell)-(1-\ell^{-\lfloor x^\delta\rfloor})(2x)^{\lfloor x^\delta\rfloor}\right|\\
\leq & \left| \# \mathcal{A}_{\leq x}(\ell,\lfloor x^\delta\rfloor ) - (1-\ell^{-\lfloor x^\delta\rfloor})(2x)^{\lfloor x^\delta\rfloor}\right|+\sum_{t=1}^{\lfloor x^\delta\rfloor-1} \# \mathcal{A}_{\leq x}(\ell, t)\\
\leq & c\lfloor x^\delta\rfloor x^{\lfloor x^\delta\rfloor-1}+ \sum_{t=1}^{\lfloor x^\delta\rfloor-1}\left( ct x^{t-1} + 2^t(1-\ell^{-t})x^t\right)\\
\leq & c\lfloor x^\delta\rfloor x^{\lfloor x^\delta\rfloor-1}+ (\lfloor x^\delta\rfloor -1)\left( c(\lfloor x^\delta\rfloor -1) x^{\lfloor x^\delta\rfloor-2}\right) + \frac{(2x)\left((2x)^{\lfloor x^\delta\rfloor-1}-1\right)}{2x-1}\\
\leq & c\lfloor x^\delta\rfloor x^{\lfloor x^\delta\rfloor-1}+  c x^{\lfloor x^\delta\rfloor} + 2(2x)^{\lfloor x^\delta\rfloor-1}.\\
\end{split}\]
After replacing $c$ with a larger absolute constant, we find that each of the three terms in the above sum are $<\frac{c}{3} (2x)^{\lfloor x^\delta\rfloor-1}$, and thus, their sum is $<c (2x)^{\lfloor x^\delta\rfloor-1}$. This proves the result.
\end{proof}
Set $g(\ell, t)$ to be the number of tuples $\bar{\alpha}=(\bar{\alpha}_1,\dots, \bar{\alpha}_t)\in \F_\ell^t$ such that 
\begin{itemize}
    \item $\bar{\alpha}\neq 0$, 
    \item $\sum_i \bar{\alpha}_i^2=0$. 
\end{itemize}

Note that by the proof of Theorem \ref{thm:mu0lambda1}, we have that
\begin{equation}\label{g ell t eqn}g(\ell, t)=\begin{cases}
\ell^{t-1}+(-1)^{\frac{t(\ell-1)}{4}}(\ell-1)\ell^{\frac{t}{2}-1}-1,& \text{ if }t \text{ is even};
\\ 
\ell^{t-1}-1, & \text{ if }t\text{ is odd.}
\end{cases}\end{equation}

\begin{lemma}\label{lemma 4.24}
With respect to the notation above, 
\[\left|\#\mathcal{T}_{\leq x}(0, 1, \ell, t)-\left(1-\frac{g(\ell, t)}{\ell^t-1}\right)(1-\ell^{-t})(2x)^t\right|<2^t\ell tx^{t-1}.\]
\end{lemma}
\begin{proof}
Let $\mathfrak{U}$ be the set of vectors $\bar{\alpha}\in \F_\ell^t$ such that 
\begin{itemize}
    \item $\bar{\alpha}\neq 0$, 
    \item $\sum_i \bar{\alpha}_i^2\neq 0$.
\end{itemize} Note that $\#\mathfrak{U}$ is equal to $\ell^t-1-g(\ell, t)$. Given $\bar{\alpha}\in \mathfrak{U}$, set $\mathcal{T}_{\leq x}(\bar{\alpha})$ to denote the subset of $\mathcal{T}_{\leq x}=\mathcal{T}_{\leq x}(0,1,\ell,t)$ consisting of tuples $\alpha$ that reduce to $\bar{\alpha}$ modulo-$\ell$. Let $\alpha_i'$ be the integer such that $0\leq \alpha_i'<\ell$ and that $\alpha_i'$ reduces to $\bar{\alpha}_i$. We have 
\[\# \mathcal{T}_{\leq x}(\bar{\alpha})=\prod_{i=1}^t\left(\lfloor (x-\alpha_i')/\ell\rfloor+ \lfloor (x+\alpha_i')/\ell\rfloor\right).\]
It then follows that 
\[\# \mathcal{T}_{\leq x}(\bar{\alpha})\leq \prod_{i=1}^t\left((x-\alpha_i')/\ell+  (x+\alpha_i')/\ell\right)\leq  \left(2x/\ell\right)^t.\]

Let $\pi_i:=\{(x-\alpha_i')/\ell\}+\{(x+\alpha_i')/\ell\}$, 
we have that 
\[\begin{split}\# \mathcal{T}_{\leq x}(\bar{\alpha})=\prod_{i=1}^{t} \left(2x/\ell-\pi_i\right)\geq (2x/\ell)^t-\left(\sum_{i=1}^t\pi_i\right) (2x/\ell)^{t-1}\end{split}>(2x/\ell)^t-\left(\frac{2^ttx^{t-1}}{\ell^{t-1}}\right).\]
This allows us to deduce that
\[\begin{split}& \left|\#\mathcal{T}_{\leq x}(0, 1, \ell, t)-\left(1-\frac{g(\ell, t)}{\ell^t-1}\right)(1-\ell^{-t})2^t x^t\right|\\
=&\left|\#\mathcal{T}_{\leq x}(0, 1, \ell, t)-\left(\ell^t-1-g(\ell, t)\right)\left(\frac{2x}{\ell}\right)^t \right|\\
= &  \left| \sum_{\bar{\alpha}\in \mathfrak{U}}\left(\# \mathcal{T}_{\leq x}(\bar{\alpha}) -\left(\frac{2x}{\ell}\right)^t\right) \right|\\
\leq & \sum_{\bar{\alpha}\in \mathfrak{U}} \left|\# \mathcal{T}_{\leq x}(\bar{\alpha}) -\left(\frac{2x}{\ell}\right)^t\right|\\
\leq & \left(\ell^t-1-g(\ell,t)\right)\left(\frac{2^ttx^{t-1}}{\ell^{t-1}}\right)\\
< & 2^t\ell tx^{t-1}.
\end{split}\]
\end{proof}
\begin{lemma}\label{lemma 4.25}
\[\left|\#\mathcal{T}_{\leq x}(0, 1, \ell, t)-(1-1/\ell)(1-\ell^{-t})(2x)^t\right|<2^t\ell t x^{t-1}+(2x)^t/\ell^{t/2}+(2x/\ell)^t.\]
\end{lemma}
\begin{proof}
We have the following estimates
\[
\begin{split}
    & \left|\#\mathcal{T}_{\leq x}(0, 1, \ell, t)-(1-1/\ell)(1-\ell^{-t})2^t x^t\right|\\
    \leq & \left|\#\mathcal{T}_{\leq x}(0, 1, \ell, t)-\left(1-\frac{g(\ell, t)}{\ell^t-1}\right)(1-\ell^{-t})2^t x^t\right|\\
    &+\left|\left(\frac{\ell^t-1-g(\ell, t)}{\ell^t}\right)2^t x^t-(1-1/\ell)2^t x^t\right|+(1-1/\ell)\ell^{-t}2^tx^t\\
    < & 2^t\ell t x^{t-1}+|\ell^{t-1}-1-g(\ell, t)|(2x/\ell)^t+(2x/\ell)^t\\
    < & 2^t\ell t x^{t-1}+(2x)^t/\ell^{t/2}+(2x/\ell)^t.
\end{split}
\]
In the above, the second inequality follows from Lemma \ref{lemma 4.24}, and the third inequality follows from \eqref{g ell t eqn}.
\end{proof}
\begin{proposition}\label{prop 4.23}With respect to notation above, we have that
\[\left|\#\mathcal{T}_{\leq x}(0, 1, \ell)-(1-1/\ell)(1-\ell^{-\lfloor x^\delta\rfloor})(2x)^{\lfloor x^\delta\rfloor}\right|<c\ell {\lfloor x^\delta\rfloor} (2x)^{{\lfloor x^\delta\rfloor}-1}\]
for some absolute constant $c>0$.
\end{proposition}
\begin{proof}
By Lemma \ref{lemma 4.25}, we have that
\[\begin{split}& \left|\#\mathcal{T}_{\leq x}(0, 1, \ell)-(1-1/\ell)(1-\ell^{-\lfloor x^\delta\rfloor})(2x)^{\lfloor x^\delta\rfloor} \right|\\
\leq &  \left|\#\mathcal{T}_{\leq x}(0, 1, \ell, \lfloor x^\delta\rfloor )-(1-1/\ell)(1-\ell^{-\lfloor x^\delta\rfloor})(2x)^{\lfloor x^\delta\rfloor}\right|\\
&+\sum_{t=1}^{\lfloor x^\delta\rfloor-1}\#\mathcal{T}_{\leq x}(0, 1, \ell, t )\\
\leq & 2^{\lfloor x^\delta\rfloor}\ell \lfloor x^\delta\rfloor x^{\lfloor x^\delta\rfloor-1}+(2x)^{\lfloor x^\delta\rfloor}/\ell^{{\lfloor x^\delta\rfloor}/2}+\left(2x/\ell\right)^{\lfloor x^\delta\rfloor}\\
&+\sum_{t=1}^{\lfloor x^\delta\rfloor-1}\left((1-1/\ell)2^t x^t+2^t\ell t x^{t-1}+(2x)^t/\ell^{t/2}\right)\\
\leq & 2^{\lfloor x^\delta\rfloor}\ell \lfloor x^\delta\rfloor x^{\lfloor x^\delta\rfloor-1}+(2x)^{\lfloor x^\delta\rfloor}/\ell^{{\lfloor x^\delta\rfloor}/2}+\left(2x/\ell\right)^{\lfloor x^\delta\rfloor}\\
&+\sum_{t=1}^{\lfloor x^\delta\rfloor-1}\left(2^t x^t+2^t\ell t x^{t-1}\right)\\
\leq & 2^{\lfloor x^\delta\rfloor}\ell \lfloor x^\delta\rfloor x^{\lfloor x^\delta\rfloor-1}+(2x)^{\lfloor x^\delta\rfloor}/\ell^{{\lfloor x^\delta\rfloor}/2}+\left(2x/\ell\right)^{\lfloor x^\delta\rfloor}\\
& +\frac{(2x)^{\lfloor x^\delta\rfloor}-1}{2x-1}+2\ell(\lfloor x^\delta\rfloor-1)^2 (2x)^{{\lfloor x^\delta\rfloor}-2}\\
\leq & c 2^{\lfloor x^\delta\rfloor}\ell \lfloor x^\delta\rfloor x^{\lfloor x^\delta\rfloor-1},
\end{split}\]
where $c>0$ is an absolute constant. Note that in the above sum, we used the estimate
\[\begin{split}&\sum_{t=1}^{\lfloor x^\delta\rfloor-1}2^t\ell t x^{t-1}\\
\leq &\left(\lfloor x^\delta\rfloor-1\right) \op{max}\{2^t\ell t x^{t-1}\mid t=1, \dots, \lfloor x^\delta\rfloor-1\}\\
\leq &\left(\lfloor x^\delta\rfloor-1\right)\left(2^{\lfloor x^\delta\rfloor-1}\ell (\lfloor x^\delta\rfloor-1) x^{\lfloor x^\delta\rfloor-2}\right).\end{split}\]

\end{proof}

\begin{proof}[Proof of Theorem \ref{section 3 main thm}]
It follows from Proposition \ref{prop 4.20} that
\[\left| \frac{\# \mathcal{A}_{\leq x}(\ell)}{(1-\ell^{-\lfloor x^\delta\rfloor})(2x)^{\lfloor x^\delta\rfloor}} - 1\right|< \frac{c(2 x)^{\lfloor x^\delta\rfloor -1}}{(1-\ell^{-\lfloor x^\delta\rfloor})(2x)^{\lfloor x^\delta\rfloor}}=c(1-\ell^{-\lfloor x^\delta\rfloor})^{-1}(2x)^{-1},\] thus,
\begin{equation}\label{boring eqn 1}\lim_{x\rightarrow \infty} \left( \frac{\# \mathcal{A}_{\leq x}(\ell)}{(2x)^{\lfloor x^\delta\rfloor}}\right)=1. \end{equation}

By Proposition \ref{prop 4.23}, 
\[\begin{split}\left|\frac{\#\mathcal{T}_{\leq x}(0, 1, \ell)}{(1-\ell^{-\lfloor x^\delta\rfloor})(2x)^{\lfloor x^{\delta}\rfloor}}-(1-1/\ell)\right|<& \frac{c \ell {\lfloor x^\delta\rfloor} (2x)^{{\lfloor x^\delta\rfloor}-1}}{(1-\ell^{-\lfloor x^\delta\rfloor})(2x)^{\lfloor x^\delta\rfloor}}\\
= & c (1-\ell^{-\lfloor x^\delta\rfloor})^{-1}\ell {\lfloor x^\delta\rfloor} (2x)^{-1}.\end{split}\]
Recall that $\delta<1$, and thus ${\lfloor x^\delta\rfloor} (2x)^{-1}$ goes to $0$ as $x$ goes to $\infty$. We thus find that
\[\lim_{x\rightarrow \infty} \left( \frac{\#\mathcal{T}_{\leq x}(0, 1, \ell)}{(2x)^{\lfloor x^{\delta}\rfloor}}\right)=1-\ell^{-1}, \] and the result follows from \eqref{boring eqn 1}.

\end{proof}

\section{Statistics for multigraphs with two vertices}\label{S:2-vertices}
\par In this section, we study statistical problems associated to multigraphs consisting of $2$ vertices $v_1$ and $v_2$. We study an analogue of Question \eqref{bouquet question 1} for such graphs. First, let us set up some notation. Let $p$ (resp. $q$) be the number of undirected loops starting and ending at $v_1$ (resp. $v_2$). Denote by $r$ the number of undirected edges joining $v_1$ and $v_2$. The Iwasawa invariants we study significantly depend on the choice of the set $S$. Let $e$ (resp. $g$) be the number of directed edges in $S$ starting at $v_1$ (resp. $v_2$) and ending at $v_2$ (resp. $v_1$). Note that $e+g=r$. We shall assume that $r\geq 1$. If $e=0$ and $r\geq 1$, then we replace $v_1$ with $v_2$ and thus without loss of generality assume that $e\geq 1$. Any graph with $2$-vertices and choice of $S$ corresponds to the $5$-tuple $(p,q,r,e,g)$ and let $X_{p,q,r}^{e,g}$ be the associated datum. It conveniences us to express the vector $\alpha$ as $\alpha=(a_1, \dots, a_p, b_1, \dots, b_e, c_1,\dots, c_g, d_1, \dots, d_q)$,
where $a_i$ is the value of $\alpha$ at the $i$-th loop in $S$ based at $v_1$. We set $b_i$ to be the value of $\alpha$ at the $i$-th directed edge in $S$ starting at $v_1$ and ending at $v_2$. Denote by $c_i$ the value of $\alpha$ at the $i$-th directed edge in $S$ starting at $v_2$ and ending at $v_1$. Finally, let $d_i$ be the value of $\alpha$ at the $i$-th loop in $S$ based at $v_2$. 
\par We write $\vec{a}$, $\vec{b}$, $\vec{c}$ and $\vec{d}$ to denote the vectors $(a_1,\dots,a_p)$, $(b_1, \dots, b_e)$, $(c_1, \dots, c_g)$ and $(d_1, \dots, d_q)$, respectively. Then, $\alpha$ may be represented as a matrix 
\[\alpha=\mtx{\vec{a}}{\vec{b}}{\vec{c}}{\vec{d}},\]where the $(i, j)$-th entry in this matrix represents the values that $\alpha$ takes at the various directed edges that start at $v_i$ and end at $v_j$. Thus, the vectors $\vec{a}$ and $\vec{d}$ represent the values taken at loops starting and ending at $v_1$ and $v_2$, respectively.

\par Recall that Theorem \ref{connectedness_cond} gives conditions for Assumption \ref{main assumption} to be satisfied. We fix a spanning tree $\mathfrak{T}$ and assume that $\alpha:S\rightarrow \Z_\ell$ satisfies the conditions of Theorem \ref{connectedness_cond}. Without loss of generality, assume that $\mathfrak{T}$ consists of a single edge joining $v_1$ to $v_2$ such that $S(\mathfrak{T})$ consists of the first directed edge joining $v_1$ to $v_2$. Then the condition that $\alpha$ is equal to $0$ on $S(\mathfrak{T})$ is equivalent to the condition that $b_1=0$. The condition that $\alpha(s)\in \Z_\ell^\times $ for some $s\in \Z_\ell^\times$ means that one of the coordinates of $\alpha$ is in $\Z_\ell^\times$. We say that $\alpha$ is admissible if these conditions on $\alpha$ are satisfied. It follows that if $\alpha$ is admissible, then Assumption \ref{main assumption} is satisfied. Set $t:=p+q+r-1$ and given $\alpha$, set $\alpha_0\in \Z_\ell^t$ denote the vector \[\alpha_0=(a_1, \dots, a_p, b_2, \dots, b_e, c_1, \dots, c_g, d_1, \dots, d_q),\] obtained by dropping the coordinate for $b_1$. When $\alpha$ is admissible, at least one of the coordinates of $\alpha_0$ is not divisible by $\ell$. 

For example, let $\ell=5$, $X$ the multigraph consisting of two vertices $v_{1}, v_{2}$, two loops at $v_{1}$, one loop at $v_{2}$ and two undirected edges between $v_{1}$ and $v_{2}$, and 
$$S = \{s_{1},s_{2},s_{3},s_{4},s_{5} \},$$ 
where $s_{1}, s_{2}$ are the loops at $v_{1}$, $s_{3}, s_{4}$ are both directed edges going from $v_{1}$ to $v_{2}$, and $s_{5}$ is a loop at $v_{2}$.  Let $\alpha:S \to \mathbb{Z}_{\ell}$ be defined via $s_{1}, s_{2}, s_{5} \mapsto 1$, $s_{3} \mapsto 0$, and $s_{4} \mapsto 2$, then we get the $\mathbb{Z}_{5}$-tower
$$
\parbox{.20\linewidth}{
    \centering
        \includegraphics[scale=0.25]{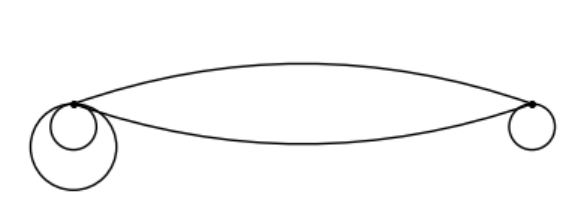}
} 
\parbox{.03\linewidth}{
    $\longleftarrow$
}
\parbox{.20\linewidth}{
    \centering
        \includegraphics[scale=0.25]{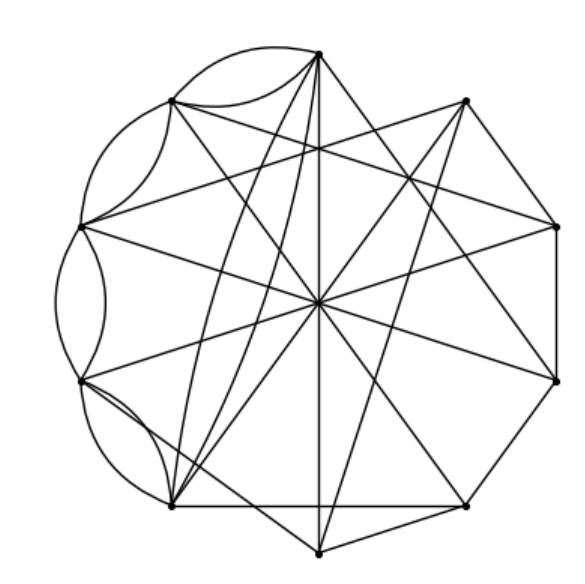}
} 
\parbox{.03\linewidth}{
    $\longleftarrow$
}
\parbox{.20\linewidth}{
    \centering
        \includegraphics[scale=0.25]{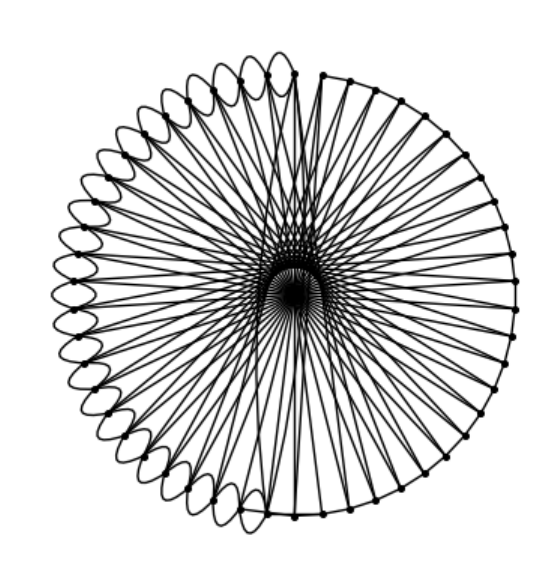}
}
\parbox{.03\linewidth}{
    $\longleftarrow$
}\hspace{0.4cm}
\parbox{.03\linewidth}{
    $\cdots$
}
$$
The power series $f(T)$ starts as follows
$$f(T) = -10T^{2} +10T^{3} -9T^{4} + \ldots $$
Thus, $\mu=0$ and $\lambda = 3$.  Using SageMath, we calculate
$$\kappa_{0} =  2, \kappa_{1} = 2 \cdot 5^{3} \cdot 31^{2}, \kappa_{2} = 2 \cdot 5^{6} \cdot 31^{2} \cdot 5351^{2} \cdot 2157401^{2}, \ldots $$
We have
$${\rm ord}_{5}(\kappa_{n}) = 3n, $$
for all $n \ge 1$.

Fixing the choice of $S$ and $\mathfrak{T}$, we parametrize the set of admissible vectors by the set of all vectors $\alpha_0\in \Z_\ell^t/\left(\ell \Z_\ell\right)^t$. We denote this set by $\mathcal{A}(\ell)$ and note that the measure equals $\ell^t-1=(1-\ell^{-t})\ell^t$. Let $\mu_{\alpha}$ and $\lambda_{\alpha}$ be the Iwasawa $\mu$ and $\lambda$-invariants of $X_{p,q,r}^{e,g}$ corresponding to the choice of voltage vector $\alpha$. Given constants $\mu_0$ and $\lambda_0$, let $\mathcal{T}(\mu=\mu_0, \lambda=\lambda_0, \ell)$ be the subset of $\mathcal{A}(\ell)$ such that $\mu_\alpha=\mu_0$ and $\lambda_\alpha=\lambda_0$. Let $\op{Prob}_\ell(\mu=\mu_0, \lambda=\lambda_0)$ denote the probability that $\mu=\mu_0$ and $\lambda=\lambda_0$ for $X_{p,q,r}^{e,g}$. More precisely, the probability is given by 
\[\op{Prob}_\ell(\mu=\mu_0, \lambda=\lambda_0):=\frac{\tau\left(\mathcal{T}(\mu=\mu_0, \lambda=\lambda_0, \ell)\right)}{\ell^t-1},\]
where $\tau$ is the Haar measure on $\Z_\ell^t$ as in \S\ref{S:bouquet}.

We are able to compute the above probability in the special case $\mu_0=0$ and $\lambda_0=1$. Since the $\lambda$-invariant is odd, these are the minimal values of the Iwasawa invariants. The main result of this section is Theorem \ref{main thm 2 vertices}. 

Recall that $D$ is the diagonal matrix given by $D=\op{diag}\left((2p+r),(2q+r)\right)$, and write $M(1+T)=\left(m_{i,j}(T)\right)_{1\leq i, j\leq 2}$ (cf. \eqref{M matrix}). We are able to explicitly calculate the four values of $m_{i,j}(T)$ as follows
\[\begin{split}
    m_{1,1}(T)=&\ r+\sum_{i=1}^p \left(2-(1+T)^{a_i}-(1+T)^{-a_i}\right)
    = r-\sum_{k\geq 2}\sum_{i=1}^p \left({a_i \choose k} +{-a_i\choose k} \right)T^k,\\
     m_{1,2}(T)= & -\left(1+\sum_{i=2}^{e}(1+T)^{b_i}+\sum_{j=1}^g(1+T)^{-c_j}\right)
     = -r-\sum_{k\geq 1} \left(\sum_{i=2}^e {b_i \choose k} +\sum_{j=1}^g {-c_j \choose k}\right)T^k,\\
m_{2,1}(T)= & -\left(1+\sum_{i=2}^{e}(1+T)^{-b_i}+\sum_{j=1}^g(1+T)^{c_j}\right)
= -r-\sum_{k\geq 1} \left(\sum_{i=2}^e {-b_i \choose k} +\sum_{j=1}^g {c_j \choose k}\right)T^k,\\
    m_{2,2}(T)=&\ r+\sum_{i=1}^q \left(2-(1+T)^{d_i}-(1+T)^{-d_i}\right)
    = r-\sum_{k\geq 2}\sum_{j=1}^q \left({d_j \choose k} +{-d_j\choose k} \right)T^k.\\
\end{split}\]

We write $f(T)=\op{det}M(1+T)=\sum_{k=0}^{\infty} \beta_k T^k$ and compute the first few values of $\beta_k$.
\begin{lemma}\label{lemma two vertices}
With respect to notation above, $\beta_n=0$ for $n\leq 1$. We find that 
\[\beta_2=\left(\sum_{i=2}^e b_i-\sum_{j=1}^g c_j\right)^2-r\left(\sum_{i=1}^p a_i^2+\sum_{j=2}^e b_j^2+\sum_{k=1}^g c_k^2+\sum_{n=1}^q d_n^2\right).\] The Iwasawa $\lambda$-invariant satisfies $\lambda\geq 1$, and the following are equivalent
\begin{enumerate}
    \item $\mu=0$ and $\lambda=1$,
    \item $\ell\nmid \beta_2$.
\end{enumerate}
\end{lemma}
\begin{proof}
In order to compute $\beta_0$, $\beta_1$ and $\beta_2$, we find the matrix $M(1+T)$ modulo $T^3$. We find that 
\[M(1+T)=\mtx{r-\xi_1 T^2}{-r-\alpha T-\xi_2 T^2}{-r+\alpha T -\xi_3 T^2}{r-\xi_4T^2}\pmod{T^3},\] where \[\begin{split}& \alpha=\sum_{i=2}^e b_i-\sum_{j=1}^g c_j,\\ 
& \xi_1= \sum_{i=1}^p a_i^2,\\ & \xi_2= \sum_{i=2}^e{b_i\choose 2} +\sum_{j=1}^g {-c_j\choose 2},\\ &\xi_3= \sum_{i=2}^e{-b_i\choose 2} +\sum_{j=1}^g {c_j\choose 2},\\ &\xi_4= \sum_{j=1}^q d_j^2.\\
\end{split}\] 
Thus, \[\begin{split}f(T)=& \op{det}M(1+T)\\
\equiv & (r-\xi_1 T^2)(r-\xi_4 T^2)-(-r+\alpha T-\xi_2 T^2)(-r-\alpha T-\xi_3 T^2)\\ \equiv & \left(\alpha^2-r\sum_{i=1}^4 \xi_i\right)T^2\mod{T^3}\\
\equiv & \left(\left(\sum_{i=2}^e b_i-\sum_{j=1}^g c_j\right)^2-r\left(\sum_{i=1}^p a_i^2+\sum_{j=2}^e b_j^2+\sum_{k=1}^g c_k^2+\sum_{n=1}^q d_n^2\right)\right)T^2\mod{T^3}.\end{split}\] In other words, $\beta_0=\beta_1=0$, and 
\[\beta_2=\left(\sum_{i=2}^e b_i-\sum_{j=1}^g c_j\right)^2-r\left(\sum_{i=1}^p a_i^2+\sum_{j=2}^e b_j^2+\sum_{k=1}^g c_k^2+\sum_{n=1}^q d_n^2\right).\]
\end{proof}
\par We study the proportion of admissible vectors $\alpha$ 
(and thus $\Z_\ell$-towers) such that $\mu=0$ and $\lambda=1$. These are the minimal values of the Iwasawa invariants. We consider the quadratic form over $\F_\ell$ in $t$ variables \[\begin{split}& Q_{p, e,g,q}\left(x_1, \dots, x_p, y_2, \dots, y_e, z_1,\dots, z_g, w_1, \dots, w_q\right) \\
:= & \left(\sum_{i=2}^e y_i-\sum_{j=1}^g z_j\right)^2-r\left(\sum_{i=1}^p x_i^2+\sum_{j=2}^e y_j^2+\sum_{k=1}^g z_k^2+\sum_{n=1}^q w_n^2\right).\\ \end{split}\]

Note that if $\ell|r$, then the quadratic form above is simply given by
\[Q_{p, e,g,q}\left(x_1, \dots, x_p, y_2, \dots, y_e, z_1,\dots, z_g, w_1, \dots, w_q\right) =\left(\sum_{i=2}^e y_i-\sum_{j=1}^g z_j\right)^2,\] and is clearly degenerate. 
\begin{proposition}\label{rank of Q}
Let $\ell$ be an odd prime number and $p,q,r,e,g$ as above. Then, the following assertions hold.
\begin{enumerate}
    \item\label{prop 5.2 1} If $\ell|r$, then $Q_{p, e,g,q}$ is a degenerate quadratic form over $\F_\ell$ of rank $1$. 
    \item\label{prop 5.2 2} If $\ell\nmid r$, then $Q_{p, e,g,q}$ is a non-degenerate quadratic form over $\F_\ell$ of rank $t$.
\end{enumerate}
\end{proposition}
\begin{proof}
Let $Q:=Q_{p,e,g,q}$ and let $\langle \cdot , \cdot \rangle$ be the associated bilinear form defined by 
\[\langle u, v \rangle = \frac{1}{2}\left( Q(u+v)-Q(u)-Q(v)\right).\]
Let \[u=\left(x_1, \dots, x_p, y_2, \dots, y_e, z_1,\dots, z_g, w_1, \dots, w_q\right)\] and \[u'=\left(x_1', \dots, x_p', y_2', \dots, y_e', z_1',\dots, z_g', w_1', \dots, w_q'\right).\] We note that
\[\begin{split}& \langle u, u' \rangle= \\
 & \left(\sum_{i=2}^e y_i-\sum_{j=1}^g z_j\right)\left(\sum_{i=2}^e y_i'-\sum_{j=1}^g z_j'\right)-r\left(\sum_{i=1}^p x_i x_i'+\sum_{j=2}^e y_j y_j'+\sum_{k=1}^g z_k z_k'+\sum_{n=1}^q w_n w_n'\right).\\ \end{split}\]

The assertion \eqref{prop 5.2 1} is clear since the associated pairing is given by 
\[\langle u, u'\rangle =\left(\sum_{i=2}^e y_i-\sum_{j=1}^g z_j\right)\left(\sum_{i=2}^e y_i'-\sum_{j=1}^g z_j'\right).\] Therefore, in order for $\langle u, \cdot \rangle$ to be identically $0$, it is necessary and sufficient for $\sum_{i=2}^e y_i-\sum_{j=1}^g z_j$ to equal $0$. The rank is therefore equal to $t-(t-1)=1$.

\par Let us prove part \eqref{prop 5.2 2}, i.e., we show that $Q$ is non-degenerate if $\ell\nmid r$.
Let $u\neq 0$ be such that $\langle u, u' \rangle=0$ for all $u'$. Thus $x_i$ and $w_n$ are all clearly $0$.  Consider $u'$ such that all coordinates are $0$ except for the coordinate $y_j'$, which is $1$. We find that 
\begin{equation}\label{u,u' 1}\langle u,u' \rangle=\left(\sum_{i=2}^e y_i-\sum_{j=1}^g z_j\right)-r y_j=0.\end{equation}
Since $\ell\nmid r$, we find that $y_j=y_k$ for $j\neq k$. Next, consider $u'$ such that all coordinates are $0$ except for the coordinate $z_j'$ which is $1$. We find that 
\begin{equation}\label{u,u' 2}\langle u,u' \rangle=-\left(\sum_{i=2}^e y_i-\sum_{j=1}^g z_j\right)-r z_j=0.\end{equation} Setting $y:=\frac{1}{r}\left(\sum_{i=2}^e y_i-\sum_{j=1}^g z_j\right)$, it follows from \eqref{u,u' 1} and \eqref{u,u' 2} that $y_i=y$, $z_j=-y$ for all $i$ and $j$.  We find that 
\[\begin{split}\langle u,u' \rangle=&\left(\sum_{i=2}^e y-\sum_{j=1}^g (-y)\right)\left(\sum_{i=2}^e  y_i'-\sum_{j=1}^g z_j'\right)-r\left(\sum_{j=2}^e y y_j' +\sum_{k=1}^g (-y)z_k' \right)\\
=&\left((r-1)y-ry\right)\left(\sum_{i=2}^e  y_i'-\sum_{j=1}^g z_j'\right)=-y\left(\sum_{i=2}^e  y_i'-\sum_{j=1}^g z_j'\right).\end{split}\]
Therefore, it follows that $y$ must be $0$ in order for the form $\langle u, \cdot \rangle$ to be identically $0$. Hence, we have shown that $u=0$, and thus, $Q$ is non-degenerate.
  
\end{proof}

\begin{theorem}\label{main thm 2 vertices}
Let $\ell$ be an odd prime and fix $(p,q,r,e,g)$ such that $e+g=r$; set $t:=p+q+r-1$. Then, the following assertions hold. 
\begin{enumerate}
    \item\label{thm 5.3 1} Suppose that $\ell\mid r$. Then, \[\op{Prob}_\ell\left(\mu=0, \lambda=1\right)=1-(1-\ell^{-t})^{-1}\left(\ell^{-1}-\ell^{-t}\right).
\]
\item\label{thm 5.3 2} Suppose that $\ell\nmid r$. Then, there is a number $\upsilon\in \{+1, -1\}$ for which
\[\op{Prob}_\ell\left(\mu=0, \lambda=1\right)=\begin{cases}1-(1-\ell^{-t})^{-1}\left(\ell^{-1}-\ell^{-t}\right),&\text{ if }t \text{ is odd;}\\
1-(1-\ell^{-t})^{-1}\left(\ell^{-1}+\upsilon (\ell-1)\ell^{-(t/2+1)}-\ell^{-t}\right),&\text{ if }t \text{ is even.}\\
\end{cases}\]
\end{enumerate}
\end{theorem}
\begin{proof}
First, we prove \eqref{thm 5.3 1}, and thus assume that $\ell \mid r$. Let $Q$ be the quadratic form defined in Proposition \ref{rank of Q}. Note that by part \eqref{prop 5.2 1} of Proposition \ref{rank of Q}, the rank of $Q$ is equal to $1$. Let $\bar{\alpha}_0$ be the reduction of $\alpha_0$ modulo $\ell$. It follows from Lemma \ref{lemma two vertices} that $\mu_{\alpha}=0$ and $\lambda_\alpha=1$ if and only if $Q(\bar{\alpha}_0)\neq 0$. Therefore, we find that 
 \[\op{Prob}_\ell\left(\mu=0, \lambda=1\right)=1-\frac{N(Q=0)-1}{\ell^t-1}.\]
 
 According to Proposition \ref{number of solutions to a quad form}, $N(Q=0)=\ell^{t-1}$. Therefore, we find that \[\op{Prob}_\ell\left(\mu=0, \lambda=1\right)=1-(1-\ell^{-t})^{-1}\left(\ell^{-1}-\ell^{-t}\right).
\]
\par For part \eqref{thm 5.3 2}, it follows from Proposition \ref{rank of Q} part \eqref{prop 5.2 2} that $Q$ is non-degenerate, and thus the rank is $t$. Setting $\upsilon$ to denote $\eta\left((-1)^{\mathfrak{t}/2}\Delta(Q)\right)$, it follows from Proposition \ref{number of solutions to a quad form} that
\[N(Q=0)=\begin{cases}  \ell^{t-1}+\upsilon(\ell-1)\ell^{t/2-1},& \text{ if }t\text{ is even;}\\
\ell^{t-1}, & \text{ if }t\text{ is odd.}\end{cases}\]
Therefore, we find that 
\[\op{Prob}_\ell\left(\mu=0, \lambda=1\right)=\begin{cases}1-(1-\ell^{-t})^{-1}\left(\ell^{-1}-\ell^{-t}\right),&\text{ if }t \text{ is odd;}\\
1-(1-\ell^{-t})^{-1}\left(\ell^{-1}+\upsilon (\ell-1)\ell^{-(t/2+1)}-\ell^{-t}\right),&\text{ if }t \text{ is even.}\\
\end{cases}\]

\end{proof}

\section{Statistics for the Iwasawa invariants of complete graphs}\label{s 6}

Given an integer $u\geq 2$, let $K_u$ be the complete multigraph with $u$ vertices, that is the multigraph with the set of vertices $V_{K_u} = \{v_1,\ldots, v_u\}$ and with exactly one directed edge $e_{i,j}$ that starts at $v_i$ and ends at $v_j$ if $i \neq j$. Furthermore, the graph has no loops. Recall that $\pi:E^+_{K_u}\to E_{K_u}$ was the quotient map taking a directed edge to its undirected class. Define $\gamma:E_{K_u} \to E^+_{K_u}$ to be the section of $\pi$ which sends the unique undirected edge joining $v_i$ to $v_j$ to $e_{i,j}$ if $i<j$ and to $e_{j,i}$ if $j<i$. Given a function $\alpha:S \to \Zl$, let us represent the values taken by $\alpha$ by a matrix $(a_{i,j})_{1\leq i,j\leq u}$ with entries defined as follows
\[a_{i,j}=\begin{cases}
\alpha(e_{i,j}) &\text{ if }i< j,\\
0 &\text{ if }i\geq j.
\end{cases}\] In this setting, the matrix $M(1+T)$ is given by $(m_{i,j}(T))_{1 \leq i,j\leq u}$ where
\[
m_{i,j}(T) = \begin{cases} u-1 & \text{if $i=j$,} \\
-(1+T)^{a_{i,j}} & \text{if $i<j$,} \\
-(1+T)^{-a_{j,i}} & \text{if $j<i$}.\end{cases}
\]

Recall that the characteristic function $f(T)$ is the determinant of $M(1+T)$. We may express $f(T)$ as a power series $f(T)=\sum_{k=2}^\infty \beta_kT^k$ by  Lemma~\ref{f(0)=0 lemma} and Theorem~\ref{lambda is odd}. The size of the matrix $M(1+T)$ increases along with $u$ making its determinant cumbersome to compute. For $u=4$ and $u=5$, we can compute explicitly
\begin{align*}
    \beta_2 =& -8\sum_{1\leq i,j \leq 4} a_{i,j}^2-8 \left( a_{1,2} a_{2,3}+ a_{1,2} a_{2,4}
    + a_{1,3} a_{3,4}+ a_{2,3} a_{3,4}\right) \\
    &+8\left( a_{1,2} a_{1,3}+a_{1,2} a_{1,4} + a_{1,3} a_{1,4}+a_{1,3} a_{2,3}+a_{2,4} a_{1,4}+a_{1,4} a_{3,4}+a_{2,4} a_{2,3}+a_{2,4} a_{3,4}\right)
\end{align*}
and
\begin{align*}
    \beta_2= &-75 \sum_{1 \leq i,j \leq 5}a_{i,j}^2 - 50 (a_{1,2} a_{2,3} + a_{1,2} a_{2,4}+a_{1,2} a_{2,5}+a_{1,3} a_{3,4}+a_{1,3} a_{3,5}+ a_{1,4} a_{4,5} \\
    &+a_{2,3} a_{3,4}+a_{2,3} a_{3,5}+a_{2,4} a_{4,5}+a_{3,4} a_{4,5}) + 50 (a_{1,2} a_{1,3}+a_{1,2} a_{1,4}+ a_{1,2} a_{1,5}
    +a_{1,3} a_{1,4}\\
    &+ a_{1,3} a_{1,5}+ a_{1,3} a_{2,3}
    +a_{1,4} a_{1,5}+ a_{1,4} a_{2,4}+ a_{1,4} a_{3,4}+ a_{1,5} a_{2,5}+ a_{3,5} a_{1,5}
    +a_{1,5} a_{4,5}\\
    &+ a_{2,3} a_{2,4}+ a_{2,3} a_{2,5}
    + a_{2,4} a_{2,5}+ a_{2,4} a_{3,4}+ a_{3,5} a_{2,5}+ a_{2,5} a_{4,5}
    + a_{3,4} a_{3,5}+ a_{3,5} a_{4,5} ),
\end{align*}
  respectively. This leads us to introduce the following concept.
For integers $1\le i\leq j\le u$ and $1\le k\leq l\le u$. We say that the tuple $(a_{i,j}, a_{k,l})$ is a \emph{linked pair} if one of the following  three conditions are satisfied
\begin{itemize}
    \item $i=k$ and $i \neq j \neq l$, \\
    \item $j=k$ and $i \neq j \neq l$, \\
    \item $j=l$ and $i \neq j \neq k$.
\end{itemize}
For example, $(a_{1,2},a_{2,4})$, $(a_{1,2},a_{1,4})$ and $(a_{1,3},a_{2,3})$ are linked pairs while $(a_{3,4},a_{2,3})$ is not. 
\begin{definition}
Let $\Pi$ be the following set
\[
\{(a_{i,j}, a_{k,l}) : \text{$1 \leq i,j,k,l \leq u$, $(a_{i,j}, a_{k,l})$ is a linked pair and $i \leq j \leq k \leq l$} \}
\]
and let $\Pi^c$ be the complement of $\Pi$ inside the following set of all linked pairs
\[
\{(a_{i,j}, a_{k,l}) : \text{$1 \leq i,j,k,l \leq u$ and $(a_{i,j}, a_{k,l})$ is a linked pair} \}.
\]
\end{definition}
Computations for values of $u$ up to $7$ suggest that the value of  $\beta_2$ may be given by the following general formula: 
\[-(u-2)u^{u-3}\sum_{1\leq i,j\leq u}a_{i,j}^2 + 2u^{u-3}\left( \sum_{\left(a_{i,j}, a_{k,l}\right) \in \Pi^c}a_{i,j}\cdot a_{k,l}-\sum_{\left(a_{i,j}, a_{k,l}\right)\in \Pi}a_{i,j}\cdot a_{k,l}\right).
\]

In order to write down an explicit formula for the first non-zero coefficient of $f(T)$, we will make simplifying assumptions. First, let us suppose the case where all the terms $a_{i,j}$ are zero except for one. This is referred to as a single voltage assignment and we note that computations in this case have been carried out in \cite{Gonet:2021}. For simplicity, let $a:=a_{12}$ be the non-zero value and assume that $\ell \nmid a$. Then, the derived multigraphs $K_u(\Z/\ell^n \mathbb{Z},S,\alpha_{/n})$ are connected for all $n$, which   follows from Theorem~\ref{connectedness_cond} by taking $\mathfrak{T}$ to be any spanning tree that does not contain an edge joining $v_1$ to $v_2$. Under these assumptions, only the coefficient $\beta_2$ will be needed to compute the distribution of the Iwasawa invariants. If one wanted to treat the case of more general voltage assignments, then one would need to compute $\beta_n$ for $n\geq 3$. 

\begin{lemma}\label{lem:char poly Kn}
Let $\alpha$ be a single voltage assignment described  above. The $\mu$-invariant of $f(T)$ is $\ord_\ell((u-2)u^{u-3})$ and its $\lambda$-invariant is $1$.
\end{lemma}

\begin{proof}
By \cite[Theorem 5]{Gonet:2021}, the determinant of $M(1+T)$ (defined as in  \eqref{M matrix}) is given by the formula
\[
-(u-2)u^{u-3}\left( (1+T)^a-1 \right)^2 (1+T)^{-a}.
\]
Since $(1+T)^{-a}$ is a unit in the Iwasawa algebra, multiplying $f(T)$ by $(1+T)^a$ does not affect its Iwasawa invariants. Thus, 
\[
\mu(f(T)) = \ord_\ell(-(u-2)u^{u-3})+2\mu((1+T)^a-1)=\ord_\ell(-(u-2)u^{u-3})
\]
because $\ell \nmid a$ by assumption. This also shows that $\lambda(f(T))=1$ since the $\ell$-adic valuation of $\beta_2$ is exactly $\ord_\ell(-(u-2)u^{u-3})$ and $(1+T)^a -1=Tu(T)$, where $u(T)$ is a unit in the Iwasawa algebra. 
\end{proof}

There is another special case where the determinant of $M(1+T)$ can be computed explicitly. Suppose that $\alpha^\prime$ is a voltage assignment with $a_{1,j}=a$ for all $2 \leq j \leq u-1$ and zero for all other $a_{j,i}$ $(0\leq i\leq j \leq u)$. The derived multigraphs are again seen to be connected for all $n$ by choosing a suitable spanning tree $\mathfrak{T}$ in Theorem \ref{connectedness_cond}. One can choose the spanning tree $\mathfrak{T}$ to pass through all vertices in the following order $v_1, v_t, v_{t-1}, \dots, v_2$. Such a tree must avoid all edges joining $v_1$ to $v_{j}$ for $2\leq j \leq t-1$.

\begin{lemma}\label{lem:char poly Kn 2}
Let $\alpha^\prime$ be the voltage assignment described above. The $\mu$-invariant of $f(T)$ is $\ord_\ell ( (u-2)u^{u-3})$ and its $\lambda$-invariant is $1$.
\end{lemma}

\begin{proof}
We set $\mathfrak{t} := (1+T)^a$ in order to ease notation. With respect to the voltage assignment $\alpha'$, we find that $M(1+T)$ is given by
\[
M(1+T) = \begin{bmatrix} u-1 & -\mathfrak{t} & -\mathfrak{t} & \cdots & -\mathfrak{t} & -1 \\ -\mathfrak{t}^{-1} & u-1 & -1 & \cdots & -1 & -1 \\
-\mathfrak{t}^{-1} & -1 & u-1 & \cdots & -1 & -1 \\
\vdots & \vdots & \vdots & \ddots & \vdots & \vdots \\
-\mathfrak{t}^{-1} & -1 & -1 & \cdots & u-1 & -1 \\
-1 & -1 & -1 & \cdots & -1 & u-1\end{bmatrix}.
\]
Define the two $1 \times u$ vectors $w$ and $v$ and the $u\times u$ matrix $C$ by 
\[
w:=\begin{bmatrix} \mathfrak{t} \\ 1 \\ \vdots \\ 1 \end{bmatrix}, v:= \begin{bmatrix} -\mathfrak{t}^{-1} \\ -1 \\ \vdots \\ -1 \end{bmatrix}, C:= \begin{bmatrix} u & 0 & \cdots& 0 & \mathfrak{t}-1 \\
0 & u & \cdots &0 & 0 \\
\vdots & \vdots & \ddots & \vdots & \vdots \\
0 & 0 & \cdots & u & 0 \\
\mathfrak{t}^{-1}-1 & 0& \cdots & 0 & u \end{bmatrix}.
\]
Then, $M(1+T) = wv^t + C$ and the matrix determinant lemma gives $\det M(1+T) = (1+v^tC^{-1}w)\det C$. To find the inverse of $C$, we directly compute the minors of $C$ and find
\[
C^{-1}= \begin{bmatrix} \frac{-u}{(\mathfrak{t}-1)(\mathfrak{t}^{-1}-1)-u^2} & 0 & \cdots& 0 & \frac{\mathfrak{t}-1}{(\mathfrak{t}-1)(\mathfrak{t}^{-1}-1)-u^2} \\
0 & \frac{1}{u} & \cdots &0 & 0 \\
\vdots & \vdots & \ddots & \vdots & \vdots \\
0 & 0 & \cdots & \frac{1}{u} & 0 \\
\frac{\mathfrak{t}^{-1}-1}{(\mathfrak{t}-1)(\mathfrak{t}^{-1}-1)-u^2} & 0& \cdots & 0 & \frac{-u}{(\mathfrak{t}-1)(\mathfrak{t}^{-1}-1)-u^2} \end{bmatrix}
\]
and $\det C = -u^{u-2} \left( (\mathfrak{t}-1)(\mathfrak{t}^{-1}-1)-u^2 \right)$. Thus,
\[
\det M(1+T) = -(u-2)u^{u-3}(\mathfrak{t}+\mathfrak{t}^{-1})+2(u-2)u^{u-3}.
\]
From this expression, we see that $\beta_0=\beta_1=0$ and 
\[
\beta_2 = -(u-2)u^{u-3}\left( \binom{a}{2}+\binom{-a}{2} \right) = -(u-2)u^{u-3}a^2.
\]
Since $\ell \nmid a$, $f(T) = -(u-2)u^{u-3}T^2 U(T)$ where $U(T)$ is a unit. Therefore, $\mu(f(T)) = \ord_\ell ((u-2)u^{u-3})$ and $\lambda(f(T))=1$ as can be seen from the fact that the $\ell$-adic valuation of $\beta_2$ is $\mu(f(T))$.
\end{proof}

For example, if we take $\ell=3$, $X = K_{4}$ and $\alpha:S \to \mathbb{Z}_{\ell}$ defined via $s_{1},s_{3}, s_{5}, s_{6} \mapsto 0$, $s_{2} \mapsto 1$, and $s_{4} \mapsto 1$, where 
$$s_{1} = e_{1,4}, s_{2} = e_{1,3}, s_{3} = e_{2,4}, s_{4} = e_{1,2}, s_{5} = e_{2,3}, s_{6} = e_{3,4},$$
then we get the $\mathbb{Z}_{3}$-tower
$$
\parbox{.20\linewidth}{
    \centering
        \includegraphics[scale=0.25]{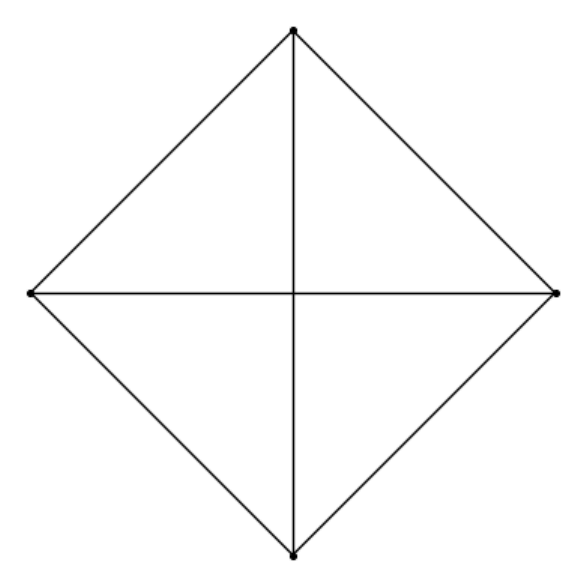}
} 
\parbox{.03\linewidth}{
    $\longleftarrow$
}
\parbox{.20\linewidth}{
    \centering
        \includegraphics[scale=0.25]{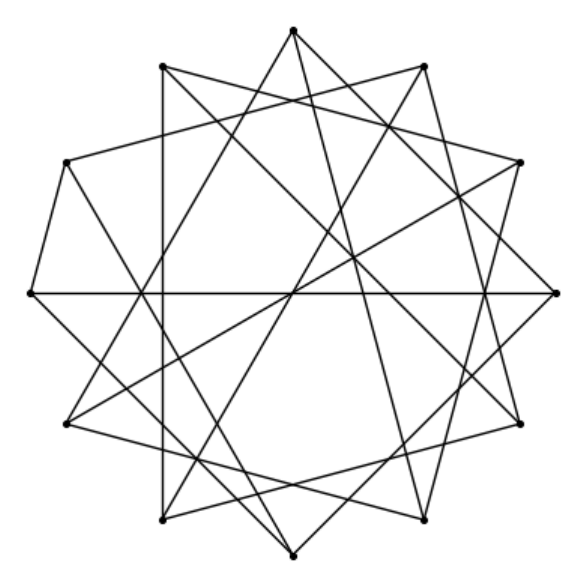}
} 
\parbox{.03\linewidth}{
    $\longleftarrow$
}
\parbox{.20\linewidth}{
    \centering
        \includegraphics[scale=0.25]{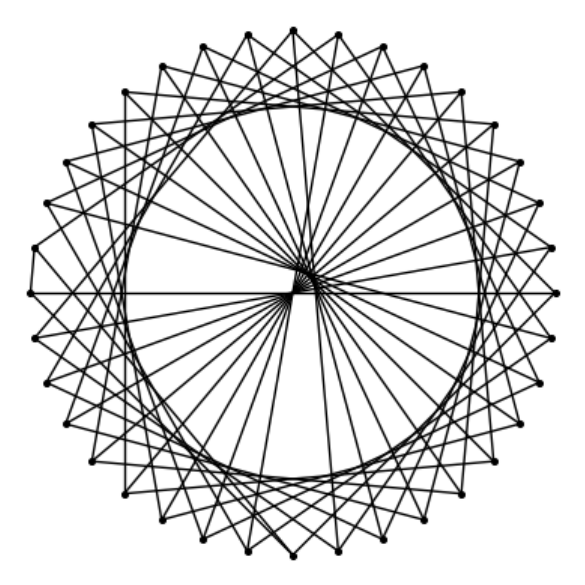}
}
\parbox{.03\linewidth}{
    $\cdots$
}
$$
The power series $f(T)$ starts as follows
$$f(T) = -8T^{2} +8T^{3} -8T^{4} + \ldots $$
Thus, $\mu=0$ and $\lambda = 1$ as expected by Lemma \ref{lem:char poly Kn 2}.  Using SageMath, we calculate
$$\kappa_{0} =  2^{4}, \kappa_{1} = 2^{10} \cdot 3, \kappa_{2} = 2^{28} \cdot 3^{2}, \kappa_{3} = 2^{82} \cdot 3^{3}, \ldots,$$
and we have
$${\rm ord}_{3}(\kappa_{n}) = n, $$
for all $n \ge 1$.

Fix $a\in \Z_\ell$ throughout and not divisible by $\ell$, set $\alpha(a,u)$ to denote the single voltage assignment with $a_{1,2}=a$ and $\alpha^\prime(a,u)$ to denote the voltage assignment with $a_{1,j}=a$ for $2\leq j \leq u-1$ and zero for the other $a_{i,j}$ with $1 \leq i \leq j \leq u$. Note that $\alpha=\alpha(a,u)$ and $\alpha^\prime = \alpha^\prime (a,u)$ depend on the choice of $a$ and $u$. We now study the variation of $\mu_{u,\alpha}=\mu\left(K_u,\alpha(a,u)\right)$ and $\lambda_{u,\alpha}=\lambda\left(K_u,\alpha(a,u)\right)$ in the case where $a\in \Z_\ell^\times$ is fixed and $u$ tends to infinity. Let $\mathcal{T}_{\leq x}(\mu_0,\lambda_0, \alpha,\ell)$ be the set of non-negative integers $u$ smaller or equal to $x$ such that $\mu_{u,\alpha}=\mu$ and $\lambda_{u,\alpha}=\lambda$. Also define these objects for $\alpha$ replaced by $\alpha^\prime$. Let $\mathcal{A}_{\leq x}$ be the set of non-negative integers $u$ smaller or equal to $x$. Note that the dependence on the choice of $a$ is suppressed in our notation.

\begin{theorem}\label{thm:complete}
Let $a\in \Z_\ell^\times$ and $\gamma \in \{\alpha,\alpha^\prime\}$. Then, we have that
\[
\lim_{x \to \infty} \frac{\#\mathcal{T}_{\leq x}(\mu,\lambda,\gamma,\ell)}{\#\mathcal{A}_{\leq x}} = \begin{cases}\frac{\ell-2}{\ell}, \quad &\text{if $\mu=0$ and $\lambda=1$};\\ \frac{\ell-1}{\ell^{\mu+1}}, &\text{if $\mu >0$ and $\lambda=1$};\\ 0, &\text{if $\lambda \neq 1$}.\end{cases}
\]
\end{theorem}

\begin{proof}
We set $\mathcal{T}_{\leq x}(\mu,\lambda,\ell)$ to denote $\mathcal{T}_{\leq x}(\mu,\lambda,\gamma,\ell)$. 
Lemma~\ref{lem:char poly Kn} and Lemma~\ref{lem:char poly Kn 2} assert that $\lambda$ is always equal to $1$ for both voltage assignments and that $\mu_{u,\gamma}=\ord_\ell(u-2)+(u-3)\ord_\ell(u)$. Suppose that $\mu_{u,\gamma}=\mu>0$ and $\lambda_{u,\gamma}=1$. Note that this is equivalent to
$$\ord_\ell(u-2)+(u-3)\ord_\ell(u)=\mu.$$ Choose $u$ large enough so that $u\geq \mu+2$. Since $\ell \neq 2$,
\[
\ord_\ell(u-2)+(u-3)\ord_\ell(u)=\mu \Leftrightarrow \ord_\ell(u-2)=\mu.
\]
Note that the number of non-negative integer with $\ell$-adic valuation equal to $\mu$ smaller or equal to $x$ is greater or equal to $\lfloor x/\ell^{\mu+1} \rfloor (\ell-1)$.
We have the following estimations:
\begin{align*}
   \#\mathcal{T}_{\leq x}(\mu,1,\ell) &= \#\mathcal{T}_{\leq \mu+1}(\mu,1,\ell)+\#\mathcal{T}_{\mu+2 \leq t \leq x}(\mu,1,\ell)\\
   &\leq \#\mathcal{T}_{\leq \mu+1}(\mu,1,\ell) + \left( \frac{x}{\ell^{\mu+1}}+1\right)(\ell-1), 
\end{align*}
\[
   \#\mathcal{T}_{\leq x}(\mu,1,\ell) \geq \#\mathcal{T}_{\mu+2 \leq t \leq x}(\mu,1,\ell) \geq \left( \frac{x}{\ell^{\mu+1}}-1 \right)(\ell-1)-\left( \frac{\mu+2}{\ell^{\mu+1}}+1\right)(\ell-1)
\]
and $x-1 \leq \#\mathcal{A}_{\leq x}\leq x$. Therefore,
\begin{align*}
\frac{\left( \frac{x}{\ell^{\mu+1}}-1 \right)(\ell-1)-\left( \frac{\mu+2}{\ell^{\mu+1}}+1\right)(\ell-1)}{x} &\leq \frac{\#\mathcal{T}_{\leq x}(\mu,1,\ell)}{\#\mathcal{A}_{\leq x}}\\
&\leq \frac{\#\mathcal{T}_{\leq \mu+1}(\mu,1,\ell) + \left( \frac{x}{\ell^{\mu+1}}+1\right)(\ell-1)}{x-1}.
\end{align*}
Letting $x \to \infty$, the result follows.

For the case $\mu_{u,\gamma}=0$ and $\lambda_{u,\gamma}=1$, we need both $\ord_\ell(u-2)$ and $\ord_\ell(u)$ to be zero. We find
\[
\frac{\left( \frac{x}{\ell}-1\right)(\ell-2)}{x} \leq \frac{\#\mathcal{T}_{\leq x}(\mu,1,\ell)}{\#\mathcal{A}_{\leq x}} \leq \frac{\left( \frac{x}{\ell}+1\right)(\ell-2)}{x-1}.
\]
We get the stated proportion of the theorem by letting $x$ goes to infinity.
\end{proof}

\begin{remark}
The sum of densities in the statement of Theorem \ref{thm:complete} for all values of $(\mu,\lambda)$ is equal to $1-1/\ell$ and not $1$. This is indeed possible for a sequence of numbers. For instance, consider the sequence $g_i=i$. Then for any value of $j>0$, the density of values $i$ for which $g_i=j$ is equal to $0$, and thus the sum of densities over all values of $j$ is $0$, not $1$.
\end{remark}

\begin{remark}
All the $\mathbb{Z}_{\ell}$-towers over complete graphs studied in this section have $\lambda = 1$, but we note in passing that there are examples where $\lambda \neq 1$.  For example, if we take $\ell=3$, $X = K_{4}$ and $\alpha:S \to \mathbb{Z}_{\ell}$ defined via $s_{4}, s_{5}, s_{6} \mapsto 0$, $s_{1} \mapsto 1$, $s_{2} \mapsto 2$, and $s_{3} \mapsto 4$, where 
$$s_{1} = e_{1,4}, s_{2} = e_{1,3}, s_{3} = e_{2,4}, s_{4} = e_{1,2}, s_{5} = e_{2,3}, s_{6} = e_{3,4},$$
then we get the $\mathbb{Z}_{3}$-tower
$$
\parbox{.20\linewidth}{
    \centering
        \includegraphics[scale=0.25]{./ex1_0.pdf}
} 
\parbox{.03\linewidth}{
    $\longleftarrow$
}
\parbox{.20\linewidth}{
    \centering
        \includegraphics[scale=0.25]{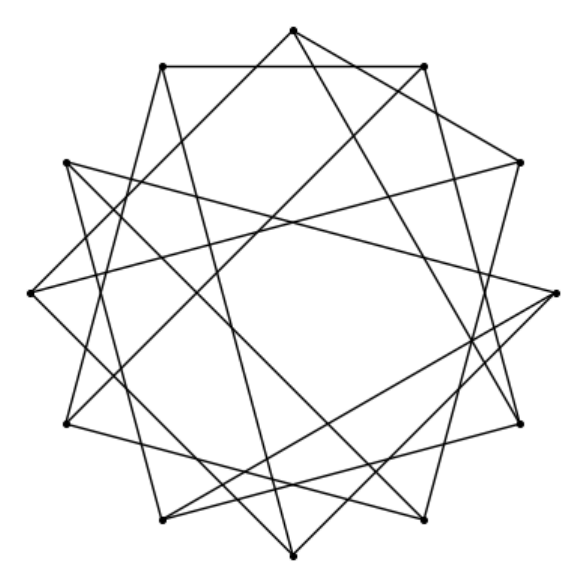}
} 
\parbox{.03\linewidth}{
    $\longleftarrow$
}
\parbox{.20\linewidth}{
    \centering
        \includegraphics[scale=0.25]{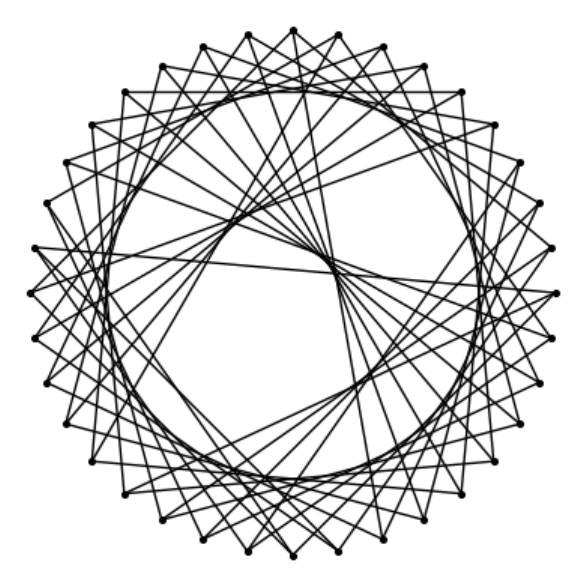}
}
\hspace{0.4cm}
\parbox{.03\linewidth}{
    $\cdots$
}
$$
The power series $f(T)$ starts as follows
$$f(T) = -120T^{2} +120T^{3} -252T^{4} + 384T^{5} - 578T^{6} + \ldots $$
Thus, $\mu=0$ and $\lambda = 5$.  Using SageMath, we calculate
$$\kappa_{0} =  2^{4}, \kappa_{1} = 2^{8} \cdot 3^{3}, \kappa_{2} = 2^{8} \cdot 3^{8} \cdot 11^{6}, \kappa_{3} = 2^{8} \cdot 3^{13} \cdot 11^{6} \cdot 13931^{2} \cdot 19996201^{2}, \ldots $$
We have
$${\rm ord}_{3}(\kappa_{n}) = 5n - 2, $$
for all $n \ge 1$.

\end{remark}

\bibliographystyle{alpha}
\bibliography{references}
\end{document}